\input amstex
\documentstyle{amsppt}
%\magnification 1200
\NoBlackBoxes

\TagsOnRight

\def\cal{\Cal}
\def\AA{{\cal A}}

\def\HH{{\cal H}}
\def\MM{{\cal M}}
\def\NN{{\cal N}}
\def\JJ{{\cal J}}
\def\UU{{\cal U}}
\def\SS{{\cal S}}
\def\KK{{\cal K}}
\def\FF{{\cal F}}
\def\LL{{\cal L}}

\def\GG{{\cal G}}
\def\Z{{\Bbb Z}}
\def\C{{\Bbb C}}
\def\R{{\Bbb R}}

\def\Q{{\Bbb Q}}
\def\e{{\epsilon}}

\def\n{\noindent}
\def\part{{\partial}}

\rightheadtext{Hamiltonian diffeomorphism group} \leftheadtext{
Yong-Geun Oh }

\topmatter
\title
Chain level Floer theory and Hofer's geometry of the Hamiltonian
diffeomorphism group
\endtitle
\author
Yong-Geun Oh\footnote{Partially supported by the NSF Grant \#
DMS-9971446 and by a grant of the Korean Young Scientist Prize
\hskip8.5cm\hfill}
\endauthor
\address
Department of Mathematics, University of Wisconsin, Madison, WI
53706, ~USA \& Korea Institute for Advanced Study, 207-43
Cheongryangri-dong Dongdaemun-gu Seoul 130-012, KOREA
\endaddress

\abstract In this paper we first apply the chain level Floer
theory to the study of Hofer's geometry of Hamiltonian
diffeomorphism group in the cases without quantum contribution: 
we prove that any quasi-autonomous Hamiltonian path on 
weakly exact symplectic manifolds or any
autonomous Hamiltonian path on arbitrary symplectic manifolds is
length minimizing in its homotopy class with fixed ends, as long
as it has a fixed maximum and a fixed minimum which are not
over-twisted and all of its contractible periodic
orbits of period less than one are sufficiently $C^1$-small. Next we give a
construction of new invariant norm of Viterbo's type on the
Hamiltonian diffeomorphism group of arbitrary compact symplectic
manifolds.
\endabstract


\endtopmatter

\document

\bigskip

\centerline{\bf Contents} \medskip

\n \S1. Introduction \smallskip

\n \S2. Normalization of Hamiltonians and the action spectrum 
\smallskip
  
\n \S3. Floer homology with real filtrations \par

 \hskip 0.1cm 3.1. Behavior of filtration under the chain homotopy
 \par

 \hskip 0.1cm 3.2. Adiabatic homotopy and adiabatic chain map
 \smallskip

\n \S4. $C^2$-small Hamiltonians and local Floer complex \par

 \hskip 0.1cm 4.1. Local Floer homology \par

 \hskip 0.1cm 4.2. $\text{Fix }\phi^1_G$ versus $\Delta \cap
 \text{graph }\phi^1_G$: comparison of two Floer homology
 \smallskip

\n \S5. Calculation
\smallskip

\n \S6. Handle sliding lemma
\smallskip

\n \S7. Non-pushing down lemma and existence 
\smallskip

\n \S8. Construction of spectral invariants
\smallskip

\n Appendix

\bigskip

\head  \bf \S1. Introduction
\endhead
In [H], Hofer introduced an invariant pseudo-norm on the group
$\HH am(M,\omega)$ of compactly supported Hamiltonian
diffeomorphisms of the symplectic manifold $(M,\omega)$ by putting
$$
\|\phi\| = \inf_{H\mapsto \phi} \|H\| \tag 1.1
$$
where $H \mapsto \phi$ means that $\phi= \phi_H^1$ is the time-one
map of Hamilton's equation 
$$
\dot x = X_H(x),
$$
and $\|H\|$ is the function defined by
$$
\|H\| = \int _0^1 \text{osc }H_t \, dt = \int_0^1(\max H_t -\min
H_t)\, dt. \tag 1.2
$$
He also proved that (1.1) is non-degenerate for the case $\C^n$
with respect to the standard symplectic structure. Subsequently,
Polterovich [Po1] and Lalonde-McDuff [LM] proved the
non-degeneracy for the case of rational symplectic manifolds and
in complete generality,  respectively. We also refer to [Ch] for
the proof in the case of tame symplectic manifolds based on the
Floer homology theory of Lagrangian intersections and its
simplification to [Oh4].

The invariant norm (1.1) induces a bi-invariant distance on $\HH
am(M,\omega)$ by
$$
d(\phi,\psi) := \|\phi\psi^{-1}\|
$$
which is the Finsler distance induced by the invariant Finsler
norm
$$
\|h\| = \max h - \min h \tag 1.3
$$
on the Lie algebra $C^{\infty}(M)/\R \simeq T_{id}\HH
am(M,\omega)$ of the group $\HH am(M,\omega)$. A natural problem
of current interest in the literature is the study of {\it
geodesics} in this Finsler manifold.

Hofer proved that the path of any {\it autonomous} Hamiltonian on
$\C^n$ is length minimizing as long as the corresponding
Hamilton's equation has no non-constant time-one periodic orbit.
This result was generalized in [MS] on general symplectic
manifolds for the case of {\it slow} autonomous Hamiltonians among
the paths homotopic with fixed ends: According to [En], [MS] and
[Mc], an autonomous Hamiltonian is called slow if it has no
non-constant contractible periodic orbit of period less than 1 and
the linearized flow at each fixed point is {\it not over-twisted}
i.e., has no closed trajectory of period less than one

We call two Hamiltonians $G$ and $F$ are called equivalent if
there exists a family $\{F^s\}_{0\leq s\leq 1}$ such that
$$
\phi^1_{F^s} = \phi^1_G 
$$
for all $s \in [0,1]$. We denote $G \sim F$ in that case and say
that two Hamiltonian paths $\phi^t_G$ and $\phi^t_F$ are homotopic
to each other with fixed ends, or just homotopic to each other
when there is no danger of confusion.

In the present paper, we study length minimizing property of the
{\it quasi-autonomous} Hamiltonian path: Such a Hamiltonian path
was proven to be geodesics in the sense of Finsler geometry [LM]
(up to time reparametrization). We refer to [Po2] for the precise
variational definition of geodesics from the first principle and
an elegant proof of this latter fact. We will just borrow theorems
from [LM] or [Po2] for a concrete description of geodesics in
terms of quasi-autonomous Hamiltonians.
\medskip

\n{\bf Definition 1.1}~ A Hamiltonian $H$ is called {\it
quasi-autonomous} if there exists two points $x_-, \, x_+ \in M$
such that
$$
H(x_-,t) = \min_x H(x,t), \quad H(x_+,t) = \max_x H(x,t)
$$
for all $t\in [0,1]$.
\medskip
It has been proven in [BP], [LM], [Po2] that a path $\{\phi^t\}$
is a geodesic in the variational sense iff the corresponding
Hamiltonian $H$ is locally quasi-autonomous. Based on this
theorem, we just say that a geodesic is the Hamiltonian path
generated by a locally quasi-autonomous Hamiltonian.

We now recall Lalonde-McDuff's necessary condition on the
stability of geodesics. In [Corollary 4.11, LM], Lalonde-McDuff
proved that for a generic $\phi$ in the sense that all its fixed
points are isolated, any stable geodesic $\phi_t, \, 0 \leq t \leq
1$ from the identity to $\phi$ must have at least two fixed points
at which the linearized isotopy has no non-constant closed
trajectory in time less than 1 in the sense of Definition 1.2
below.

\medskip

\n{\bf Definition 1.2.}~ Let $H: M \times [0,1] \to \R$ be a
Hamiltonian which is not necessarily time-periodic and $\phi_H^t$
be its Hamiltonian flow. \par

\roster \item We call a point $p\in M$ a {\it time $T$ periodic
point} if $\phi_H^T(p)=p$. We call $t \in [0,T] \mapsto
\phi_H^t(p)$ a {\it contractible time $T$-periodic orbit} if it is
contractible. \par

\item When $H$ has a fixed critical point $p$ over $t \in
[0,T]$, we call $p$ {\it over-twisted} as a time $T$-periodic
orbit if its linearized flow $d\phi_H^t(p); \, t\in [0,T]$ on
$T_pM$ has a closed trajectory of period less than $1$.
\endroster
\medskip

The following is the main result of the present paper.

\proclaim{Theorem I} Suppose that $G$ is a
quasi-autonomous Hamiltonian such that 
\smallskip

(i) all contractible periodic orbit of period less than one 
are sufficiently $C^1$-small,
\par

(ii) it has a fixed minimum and a fixed maximum which are not
over-twisted.
\smallskip

Then its Hamiltonian path $\phi_G^t$ is length minimizing in its
homotopy class with fixed ends for $ 0 \leq t \leq 1$, in cases
\roster

\item $(M,\omega)$ is weakly exact, i.e., $\omega|_{\pi_2(M)} = 0$
or \par

\item $G$ is autonomous.
\endroster
\endproclaim

\medskip

The case (1) extends the result by Siburg [Si] on $\R^{2n}$, and
(2) extends Entov's [En] and
Lalonde-McDuff-Slimowitz's result [MS] for the {\it slow}
autonomous case in several ways: first, it removes the slowness
assumption in the case of autonomous Hamiltonian. Secondly it
allows both time-dependent Hamiltonians and appearance of 
non-constant periodic orbits. Whether
Theorem I holds in general cases, when there exists quantum
contribution, is still to be seen. 

Our proof of Theorem I will be based on the Floer homology theory
on general symplectic manifolds, which has been established by now
in the general context [FOn], [LT], [Ru]. The idea of studying
length minimizing property using the Floer theory was introduced
by Polterovich [Po2] for the case of {\it small autonomous}
Hamiltonians when the action functional is single valued as in the
case of exact symplectic manifolds. We generalize his scheme to
the case of quasi-autonomous Hamiltonian paths when the action
functional is not single valued.

We first summarize Polterovich's scheme of the proof for the case
of {\it small autonomous} Hamiltonian when the symplectic form
$\omega$ is exact, say $\omega = - d\theta$. A crucial idea behind
his scheme is to relate the norm $\|h\|= h(x^+) - h(x^-)$ with two
homologically essential critical values of the action functional
$$
\AA_h(\gamma) = \int_\gamma\theta - \int_0^1 h(\gamma(t))\, dt
$$
corresponding to the maximum and minimum points $x^+$ and $x^-$ of
the function $h$, which is precisely $-h(x^-)$ and $-h(x^+)$
respectively. This is carried out first by proving some existence
result for the Floer's continuity equation
$$
\cases \frac{\part u}{\part \tau} + J (\frac{\part u}{\part t}
- X_{L^{\rho(\tau)}}(u)) = 0 \\
u(-\infty) \in \text{Crit }k, \, u(\infty) = x^+ \\
\endcases
\tag 1.4
$$
where $L^s$ is the linear homotopy
$$
L^s = (1-s) k + s F, \,\, s \in [0,1] \tag 1.5
$$
for small autonomous Hamiltonian $k$ and arbitrary Hamiltonian $F$
with $F \sim h$, and then by making some calculations involving
the action functional and the solution of (1.4). (Similar
calculations of this sort were previously employed  by Chekanov
[Ch] and by the present author [Oh3,5].) For the existence result,
Polterovich exploits the fact that when $h$ is sufficiently small,
then the Floer complex is diffeomorphic to the Morse complex of
$h$ and so the maximum point on the compact manifold $M$ is
homologically essential, which in turn is translated into the
existence of a solution of (1.4), via the fact that the Floer
complexes of $h$ and $F$ are conjugate to each other (see
Proposition 5.3), when $F \sim h$.

When we try to use
Floer homology theory in the study of quasi-autonomous Hamiltonian
paths, the first obvious point we need to take care of is that the 
Hamiltonian may not be one-periodic . This can
be taken care of using canonical modification of Hamiltonians
into time periodic ones without changing their time-one maps and
quasi-autonomous property (see Lemma 5.2 for the precise
statements).

There are many difficulties to overcome for the non-autonomous
Hamiltonians especially when the action functional is not single-valued.
However using the full power of Floer homology theory developed by
now (in the level of chain, though) and an idea of mini-max theory
via the Floer homology developed by the author in [Oh3,5], we
again reduce the proof of Theorem I to a similar existence result
(Proposition 5.3) for (1.4) where $h$ is replaced by a 
quasi-autonomous Hamiltonian. Unlike the small autonomous case, 
such an existence result is highly non-trivial (even in the autonomous 
case) for large Hamiltonians. In fact, the method we employ to prove
the existence theorem heavily relies on the extensive chain level
Floer theory. The latter turns out to carry applicability much wider 
than as we use in the present paper and leads us to 
the construction of spectral invariants on arbitrary compact symplectic
manifolds (See \S 8 and [Oh7]).

The  proof of Theorem I will then be carried out by a continuation argument 
over the homotopy
$$
\e k \mapsto \e_0G^{\e_0} \mapsto G \mapsto F
$$
combined with a delicate mini-max argument via Floer
homology over the adiabatic homotopy.
One important point that we are exploiting in the first step
is that when the Hamiltonian is $C^2$-small as in the case of 
$\e\, G^\e$ for $\e$
sufficiently small, the Floer boundary operator is decomposed into
$$
\partial = \partial_0 + \partial^\prime
$$
where $\partial_0$ is the classical contribution and
$\partial^\prime$ is the quantum contribution (see \S5, and [Oh2]
in the context of Lagrangian intersections ). This enables us to
define the concept of local Floer homology  which is invariant under local
continuation (see [Oh1] in the context of
Lagrangian submanifolds). In general
$\partial^\prime$ is not zero, but is so {\it either when $(M,\omega)$ is
weakly exact, or when the Hamiltonian is $C^2$-small 
and autonomous} which is due to the
extra $S^1$ symmetry (see [Fl2], [FHS], [FOn], [LT]). This is one place
where we used the hypotheses in Theorem I.

The second ingredient we use in this paper is several versions of
the {\it Non pushing down lemma} culminating in Proposition 7.14. 
In fact this kind of non-pushing down lemma is the 
heart of the matter in the chain level Floer theory (see [Oh7] for more
such arguments in general).
The proofs of these Non-pushing down lemmas 
use the above hypothesis in a more serious way and also use
the concept of {\it adiabatic homotopy} and {\it adiabadic chain map}.
The third ingredient is a Floer theoretic version of
the  {\it Handle sliding lemma} (Proposition 6.3). 
These tools enable us to develop a mini-max theory of 
the action functional in the non-exact case. In the much 
simpler setting of the (weakly-)exact case where
the action functional is single valued, similar
mini-max idea was previously developed by the author in [Oh3,5]
for the Lagrangian submanifolds on the cotangent bundle, and
subsequently by Schwarz [Sc] for the Hamiltonian diffeomorphisms
on symplectically aspherical symplectic manifolds. As an application
of this mini-max theory, we prove the following construction of 
spectral invariants

\proclaim{Theorem II} For each cohomology class $0 \neq a \in H^*(M,\Q)$
and Hamiltonian $H$, there exists an invariant $\rho(H;a)$ such that
$\rho(H;a) \in \text{Spec }H$ and the assignment $H \mapsto \rho(H;a)$
is $C^0$-continuous.
\endproclaim

In a sequel [Oh7] to the present paper, we have further developed 
the techniques used here and applied them to extend the definition of
these spectral invariants to the arbitrary quantum cohomology classes
$a \in QH^*(M)$. These are then applied to give a construction of
invariant norm and to obtain a new lower bound for the Hofer norm and
to the study of length minimizing property of Hofer's geodesics. 
\medskip

We would like to thank L. Polterovich for introducing us to the
idea of studying length minimizing property of geodesics in terms
of the Floer theory during his visit of KIAS Seoul, Korea, in
April 2000 and giving us a copy of his book [Po2] before its
publication. We also thank D. McDuff for sending us the preprints
[MS] and [Mc] and informing us that the proof in [LM] already
proves local length minimizing property of geodesics once
construction of Gromov-Witten invariants on general symplectic
manifolds is established. We would also like to thank her for
several helpful e-mail communications.

\medskip

\head 
\bf \S 2. Normalization of Hamiltonians and the action spectrum 
\endhead

Let $\Omega_0(M)$ be the set of contractible loops and
$\widetilde\Omega_0(M)$ be its standard covering space in the
Floer theory. We recall the definition of this covering space from
[HS] here. Note that the universal covering space of $\Omega_0(M)$
can be described as the set of equivalence classes of the pair
$(\gamma, w)$ where $\gamma \in \Omega_0(M)$ and $w$ is a map from
the unit disc $D=D^2$ to $M$ such that $w|_{\partial D} = \gamma$:
the equivalence relation to be used is that $[\overline w \#
w^\prime]$ is zero in $\pi_2(M)$.

Following Seidel [Se], we say that $(\gamma,w)$ is {\it
$\Gamma$-equivalent} to $(\gamma,w^\prime)$ iff
$$
\omega([w'\# \overline w]) = 0 \quad \text{and }\, c_1([w\#
\overline w]) = 0 \tag 2.1
$$
where $\overline w$ is the map with opposite orientation on the
domain and $w'\# \overline w$ is the obvious glued sphere. And
$c_1$ denotes the first Chern class of $(M,\omega)$. We denote by
$[\gamma,w]$ the $\Gamma$-equivalence class of $(\gamma,w)$ and by
$\pi: \widetilde \Omega_0(M) \to \Omega_0(M)$ the canonical
projection. We also call $\widetilde \Omega_0(M)$ the
$\Gamma$-covering space of $\Omega_0(M)$. The action functional
$\AA_0: \widetilde \Omega_0(M) \to \R$ is defined by
$$
\AA_0([\gamma,w]) = -\int w^*\omega. \tag 2.2
$$
Two $\Gamma$-equivalent pairs $(\gamma,w)$ and $(\gamma,w^\prime)$
have the same action and so the action is well-defined on
$\widetilde\Omega_0(M)$. When a periodic Hamiltonian $H:M \times
(\R/\Z) \to \R$ is given, we consider the functional $\AA_H:
\widetilde \Omega(M) \to \R$ by
$$
\AA_H([\gamma,w])= \AA_0(\gamma,w) - \int H(\gamma(t),t)dt
$$
Here the sign convention is chosen to be consistent with that of
[Oh3,5],
$$
\AA_H(\gamma) = \int _\gamma \theta - \int_0^1 H(\gamma(t),t)\, dt
$$
where $\omega = -d\theta$ for the canonical one form $\theta= pdq$
on the cotangent bundle which in turn is precisely the classical
mechanics Lagrangian on the cotangent bundle.

We would like to note that {\it under this convention the maximum
and minimum are reversed when we compare the action functional
$\AA_G$ and the (quasi-autonomous) Hamiltonian $G$}.

We denote by $\text{Per}(H)$ the set of periodic orbits of $X_H$.
\medskip

\n{\bf Definition 2.1 \, [Action Spectrum].}  We define the {\it
action spectrum} of $H$, denoted as $\hbox{\rm Spec}(H) \subset
\R$, by
$$
\hbox{\rm Spec}(H) := \{\AA_H(z,w)\in \R ~|~ [z,w] \in
\widetilde\Omega_0(M), z\in \text {Per}(H) \},
$$
i.e., the set of critical values of $\AA_H: \widetilde\Omega(M)
\to \R$. For each given $z \in \text {Per}(H)$, we denote
$$
\hbox{\rm Spec}(H;z) = \{\AA_H(z,w)\in \R ~|~ (z,w) \in
\pi^{-1}(z) \}.
$$

Note that $\text {Spec}(H;z)$ is a principal homogeneous space
modeled by the period group of $(M,\omega)$
$$
\Gamma_\omega = \Gamma(M,\omega) := \{ \omega(A)~|~ A \in \pi_2(M)
\}
$$
and
$$
\hbox{\rm Spec}(H)= \cup_{z \in \text {Per}(H)}\text {Spec}(H;z).
$$
Recall that $\Gamma_\omega$ is either a discrete or a countable
dense subset of $\R$.

\proclaim\nofrills{Lemma 2.2. }~ $\hbox{\rm Spec}(H)$ is a measure
zero subset of $\R$.
\endproclaim
\demo{Proof} We first note that $\text{Spec}(H;z) \subset \R$ is a
countable subset of $\R$ for each $z$. We consider the Poincar\'e
return map in a tubular neighborhood of each $z \in
\text{Per}(H)$. More precisely, we choose a small neighborhood $V
\subset M$ of $z(0)$. We identify $V$ with $2n$-ball
$B^{2n}(\delta)$ with the  point $z(0)$ identified with the center
of the ball. Choose another ball neighborhood $V^\prime=
B^{2n}(\delta')$ with $\overline V \subset V^\prime$ such that the
(first) Poincar\'e return map denoted by
$$
R_z: V \to V^\prime; \, p \mapsto \phi_H^1(p)
$$
is well-defined.  We now define a continuous map from $V$ to the
space of piecewise smooth maps from  $S^1 \cong \R/\Z$ on $M$ as
follows: for each $p \in V$, we first follow the flow of $X_H$ and
then follow from $R_z(p)$ to $p$ by the straight line under the
identification of $V^\prime$ with $B^{2n}(\delta')$. We
reparameterize the domain of the loop by re-scaling it to be
$[0,1]$.

We denote by $z_p$ the loop corresponding to $p\in V$ constructed
as above, and by $V_z \subset \Omega_0(M)$ the image of the
assignment $p \mapsto z_p$. Obviously $z_p$ is homotopic to $z$
and so any given disc $w$ bounding $z$ can be naturally continued
to bound the loop $z_p$. We denote by $w_p$ the disc continued
from $w$ and corresponding to $p \in V$. It can be easily checked
that the function
$$
h: \pi^{-1}(V_z) \to \R; \quad h([z_p,w_p]):= \AA_H([z_p,w_p])
$$
defines a smooth function on $\pi^{-1}(V_z)$ and its critical
values comprise those of $\AA_H$ near $\text{Spec}(H;z)$. This can
be proven by writing $\AA_H([z_p,w_p])$ explicitly and by a simple
local calculation. Noting that $\pi^{-1}(V_z)$ is a finite
dimensional (in fact, $2n$ dimensional) manifold, Sard's
theorem implies that the set of critical values is a measure zero
subset in $\R$. Since a finite number of such tubular
neighborhoods together with their complement cover $M$,
$\text{Spec}(H)\subset \R$ is a finite union of measure zero
subset of $\R$ and so itself has measure zero. \qed
\enddemo

For given $\phi \in {\cal  H }am(M,\omega)$, we denote by $H
\mapsto \phi$ if $\phi^1_H = \phi$, and denote
$$
\HH(\phi) = \{ H ~|~ H \mapsto \phi \}.
$$
We say that two Hamiltonians $H$ and $K$ are equivalent if they
are connected by one parameter family of Hamiltonians
$\{F^s\}_{0\leq s\leq 1}$ such that $F^s \mapsto \phi$ i.e.,
$$
\phi_{F^s}^1 = \phi \tag 2.3
$$
for all $s
\in [0,1]$. We denote by $[H]$ the equivalence class of $H$. Then
the universal covering space $\widetilde{{\cal  H }am}(M,\omega)$
of ${\cal  H }am(M,\omega)$ is realized by the set of such
equivalence classes. 

Let $F,G \mapsto \phi$ and denote
$$
f_t = \phi_F^t, \, g_t = \phi_G^t,\, \text{and }\,  h_t = f_t
\circ g_t^{-1}. 
$$
Note that $h= \{h_t\}$ defines a loop based at the identity.
Suppose $F\sim G$ so there exists a family
$\{F^s\}_{0\leq s \leq 1} \subset \HH(\phi)$ with $F_1 =F$ and
$F_0 = G$ and satisfying (2.3). In particular $h$ defines a
contractible loop. If we denote $f^s_t = \phi_{F^s}^t$,
this family provides a natural contraction of 
the loop $h$ to the identity through
$$
\widetilde h: s \mapsto f^s\circ g^{-1}; \quad f^s\circ g^{-1}(t)
: = f^s_t\circ g_t^{-1}.
$$
which in turn provides a natural lifting of the action of the
loop $h$ on $\Omega_0(M)$ to $\widetilde\Omega_0(M)$ which we
define
$$
\widetilde h \cdot [\gamma,w] = [h\gamma, \widetilde h w] \tag 2.4
$$
where $\widetilde h w$ is the natural map from $D^2$ obtained from
identifying $\widetilde h: [0,1]\times [0,1]\to {\cal H}am
(M,\omega)$ as a map from $D^2$. 

Even when $F\sim G$ and so $h$ is not contractible,
note that the (based) loop group $\Omega({\cal  H
}am(M,\omega),id)$ naturally acts on the loop space $\Omega(M)$ by
$$
(h\cdot \gamma) (t) = h(t)(\gamma(t))
$$
where $h \in \Omega({\cal H}am (M,\omega))$ and $\gamma \in
\Omega(M)$. An interesting consequence of Arnold's conjecture is
that this action maps the particular component $\Omega_0(M) 
\subset \Omega(M)$ to itself (see e.g., [Lemma
2.2, Se]). Seidel [Lemma 2.4, Se] proves that this action (by a
based loop) can be lifted to $\widetilde\Omega_0(M)$. In this
paper, we will consider only the action by contractible loops in
$\HH am(M,\omega)$. 

We now study behavior of the action spectrum $\AA_H$ when $H$
varies. In particular, we would like to study continuity property
of certain critical values which are relevant to the uniform
minimum point of the given quasi-autonomous Hamiltonian. For
this purpose, we need to normalize the spectrum
$\text{Spec}H$. We will achieve this by restricting ourselves
to $\HH_0(\phi)$ the set of normalized 
Hamiltonians with $H \mapsto \phi$ by $\int_M H_t \, d\mu = 0$
as in [Sc]. The following is proved in [Oh6] (see [Sc]
for the symplectically aspherical case where the action fuctional
is single-valued. In this case Schwarz [Sc] proved that 
the normalization works on $\HH am(M,\omega)$ not just 
on $\widetilde{\HH am} (M,\omega)$ as long as $F, \, G \mapsto \phi$
without assuming $F\sim G$).

\proclaim{Proposition 2.3 [Theorem I, Oh6]}
Let $F,\, G \in \HH_0(\phi)$ and 
$\FF = \{F^s\}_{s\in [0,1]}$ be a path in $\HH_0(\phi)$ such
that $F^0 =G$ and $F^1 = F$. Denote $h^s_t =f^s_t \circ g_t^{-1}$ and
$h^s\cdot [z,w] = [h^s\cdot z, \widetilde h^s \cdot w]$ for a $z \in 
\text{Per}(G)$. Then the function $\chi: [0,1] \to \R$ defined by
$$ 
\chi(s) = A_{F^s}(\widetilde h^s\cdot [z,w]) 
$$
is constant. In particular, we have
$$
\text{\rm Spec}(G) = \text{\rm Spec}(F).
$$
\endproclaim
  
From now on, we will always assume that the Hamiltonian functions
are normalized so that
$$
\int_M H_t\, d\mu = 0. \tag 2.5
$$

\head\bf \S 3. Floer homology with real filtration
\endhead

\n{\it 1. Behavior of filtration under the chain map} \smallskip

For each given generic $H: M \times S^1 \to \R $, we consider the
free $\Q$ vector space over
$$
\text{Crit}\AA_H = \{[z,w]\in \widetilde\Omega_0(M) ~|~ z \in
\text{Per}(H)\}. \tag 3.1
$$
To be able to define the Floer boundary operator correctly, we
need to complete this vector space downward with respect to the
real filtration provided by the action $\AA_H([z,w])$ of the
element $[z, w]$ of (3.1). More precisely,
\medskip

\n {\bf Definition 3.1.} We call the formal sum
$$
\beta = \sum _{[z, w] \in \text{Crit}\AA_H} a_{[z, w]} [z,w], \,
a_{[z,w]} \in \Q \tag 3.2 $$ a {\it Novikov chain} if there are
only finitely many non-zero terms in the expression (3.2) above
any given level of the action. We denote by $\widetilde {CF} (H)$
the set of Novikov chains.
\medskip

Here,  we put `tilde' over $CF$ to distinguish this $\Q$ vector
space with more standard Floer complex module over the Novikov
ring in the literature. Note that this is an infinite dimensional
$\Q$-vector space in general, unless $\pi_2(M) = 0$. It appears
that for the purpose of studying Hofer's geometry this set-up of
Floer homology with real filtration on the $\Gamma$-covering
space $\widetilde \Omega_0(M)$ suits better than the more standard
Floer homology on $\Omega_0(M)$ with the Novikov ring as its
coefficient, although they provide equivalent descriptions.

Since, for the study of action changes under the chain maps, we
will frequently use the chain level property of various operators
in the Floer theory, we briefly review construction of basic
operators in the Floer homology theory [Fl2]. Let $J =
\{J_t\}_{0\leq t \leq 1}$ be a periodic family of compatible
almost complex structure on $(M,\omega)$.

For each given pair $(J, H)$, we define the boundary operator
$$
\part: \widetilde{CF}(H) \to \widetilde{CF}(H)
$$
considering the perturbed Cauchy-Riemann equation
$$
\cases
\frac{\part u}{\part \tau} + J\Big(\frac{\part u}{\part t} 
- X_H(u)\Big) = 0\\
\lim_{\tau \to -\infty}u(\tau) = z^-,  \lim_{\tau \to
\infty}u(\tau) = z^+ \\
\endcases
\tag 3.3
$$
This equation, when lifted to $\widetilde \Omega_0(M)$, defines
nothing but the {\it negative} gradient flow of $\AA_H$ with
respect to the $L^2$-metric on $\widetilde \Omega_0(M)$ induced by
the metrics $g_{J_t}: = \omega(\cdot, J_t\cdot)$ . For each given
$[z^-,w^-]$ and $[z^+,w^+]$, we define the moduli space
$\MM_J([z^-,w^-],[z^+,w^+])$ of solutions $u$ of (3.3) satisfying
$$
w^-\# u \sim w^+ \tag 3.4
$$
$\part$ has degree $-1$ and satisfies $\part\circ \part = 0$.

When we are given a family $(j, \HH)$ with $\HH = \{H^s\}_{0\leq s
\leq 1}$ and $j = \{J^s\}_{0\leq s \leq 1}$, the chain
homomorphism
$$
h_{(j,\HH)}: \widetilde{CF}(J^0,H^0) \to \widetilde{CF} (J^1,H^1)
$$
is defined by the non-autonomous equation
$$
\cases \frac{\part u}{\part \tau} +
J^{\rho_1(\tau)}\Big(\frac{\part u}{\part t}
- X_{H^{\rho_2(\tau)}}(u)\Big) = 0\\
\lim_{\tau \to -\infty}u(\tau) = z^-,  \lim_{\tau \to
\infty}u(\tau) = z^+ .
\endcases
\tag 3.5
$$
where $\rho_i, \, i= 1,2$ is functions of the type $\rho :\R \to
[0,1]$,
$$
\align
\rho(\tau) & = \cases 0 \, \quad \text {for $\tau \leq -R$}\\
                    1 \, \quad \text {for $\tau \geq R$}
                    \endcases \\
\rho^\prime(\tau) & \geq 0
\endalign
$$
for some $R > 0$.  $h_{(j,\HH)}$ has degree 0 and satisfies
$$
\part_{(J^1,H^1)} \circ h_{(j,\HH)} = h_{(j,\HH)} \circ
\part_{(J^0,H^0)}.
$$

Finally when we are given a homotopy $(\overline j, \overline
\HH)$ of homotopies with $\overline j = \{j_\kappa\}$,
$\overline\HH = \{\HH_\kappa\}$, consideration of the
parameterized version of (3.5) for $ 0 \leq \kappa \leq 1$ defines
the chain homotopy map
$$
\widetilde H : \widetilde{CF}(J^0,H^0) \to \widetilde{CF}(J^1,H^1)
$$
which has degree $+1$ and satisfies
$$
h_{(j_1, \HH_1)} - h_{(j_0,\HH_0)} = \part_{(J^1,H^1)} \circ
\widetilde H + \widetilde H \circ \part_{(J^0,H^0)}.
$$
By now, construction of these maps using these moduli spaces has
been completed with rational coefficients (See [FOn], [LT] and
[Ru]). We will freely use this advanced machinery throughout the
paper. However the main stream of the proof can be read
independently of these papers once it is understood that the
bubbling of spheres is a codimension two phenomenon, which is
exactly what the advanced machinery establishes. Therefore we do
not explicitly mention these technicalities in this paper, unless
it is absolutely necessary.

The following upper estimate of the action change can be proven by
the same argument as that of [Oh3]. Because this will be a crucial
ingredient in our proof, we include its proof here for reader's
convenience.

\proclaim\nofrills{Proposition 3.2 [Theorem 7.2, Oh3].}~ When
there are two Hamiltonians $H$ and $K$, the canonical chain map

$$
h_{HK}^{lin}: \widetilde{CF}(J,H) \to \widetilde{CF}(J,K)
$$
{\it provided by the linear homotopy} $H_s = (1-s)H + sK$ respects
the filtration
$$
h_{HK}^{lin}: \widetilde{CF}^{(-\infty,a]}(J,H) \to
\widetilde{CF}^{(-\infty,a - \int \min (K-H) dt]}(J,K) \tag 3.6
$$
and so induces the homomorphism
$$
h_{HK}^{lin}: \widetilde{HF}^{(-\infty,a]}(J,H) \to
\widetilde{HF}^{(-\infty,a - \int \min (K-H) dt]}(J,K)
$$
\endproclaim
\demo{Proof} We fix $J$ here. Let $[z^+,w^+] \in \widetilde
{CF}(K)$ and $[z^-,w^-]\in \widetilde {CF}(H)$ be given. As argued
in [Oh3],  {\it for any given solution} $u$ of (3.5) and (3.4), we
compute
$$
\AA_K([z^+,w^+])- \AA_H([z^-,w^-]) = \int^\infty_{-\infty}
\frac{d}{d\tau} (\AA_{H^{\rho(\tau)}}(u(\tau))\,d\tau.
$$
Here we have
$$
\frac{d}{d\tau} \Big(\AA_{H^{\rho(\tau)}}(u(\tau)\Big)=
d\AA_{H^{\rho(\tau)}}(u(\tau))\Big( \frac{\part u}{\part
\tau}\Big) - \int^1_0\Big(\frac{\part H^{\rho(\tau)}}{\part
\tau}\Big)(u,t)\, dt.
$$
However since $u$ satisfies (3.5), we have
$$
\align d\AA_{H^{\rho(\tau)}}(u(\tau))\Big(\frac{\part u}{\part
\tau}\Big) & = \int_0^1 \omega \Big(\frac{\part u}{\part t}
- X_{H^{\rho(\tau)}}(u), \frac{\part u}{\part \tau}\Big) \, dt \\
& = - \int_0^1\Big |\frac{\part u}{\part t} -
X_{H^{\rho(\tau)}}(u)\Big |_J \leq 0
\endalign
$$
and
$$
\align \int^1_0\Big(\frac{\part H^{\rho(\tau)}}{\part
\tau}\Big)(u,t)\, dt
& = - \int_0^1 \rho^\prime(\tau)(K - H)(u,t)dt \\
& \leq - \rho^\prime(\tau) \int_0^1 \min_x(K_t-H_t)\, dt.
\endalign
$$
Combining these and using that $\rho^\prime(\tau) \geq 0$, we have
$$
\align \AA_K([z^+,w^+]) - \AA_H([z^-,w^-])
& \leq \int_{-\infty}^\infty - \rho^\prime(\tau) \int_0^1 
\min_x(K_t-H_t)\, dtd\tau\\
& = \int_{-\infty}^\infty - \rho^\prime(\tau)\, 
d\tau \int_0^1 - \min_x(K_t-H_t)\, dt\\
& = \int_0^1 - \min_x(K_t-H_t)\, dt.
\endalign
$$
By definition of the chain map $h^{lin}_{HK}$, this finishes the
proof. \qed\enddemo

\proclaim\nofrills{Proposition 3.3 [Lemma 4.3, Oh3]. }~ For a
fixed $H$ and for a given one parameter family $\overline{J} =
\{J^s\}_{s \in [0,1]}$, the natural chain map
$$
h_{\overline J}: \widetilde{CF}(J^0,H) \to \widetilde{CF}(J^1,H)
$$
respects the filtration.
\endproclaim
\demo{Proof} A similar computation, this time using (3.2) and
(3.3) with $H$ fixed, leads to
$$
\AA_H([z^+,w^+]) - \AA_H([z^-,w^-]) =
 - \int_{-\infty}^{\infty} \int_0^1 \Big|\frac{\part u}{\part t}
  - X_{H^{\rho(\tau)}}(u)\Big|_{J^{\rho(\tau)}} \leq 0.
$$
We refer to the proof of [Lemma 4.3, Oh3] for complete details.
\qed \enddemo

We would like to remark that there is also some upper estimate for
chain maps over general homotopy or for the chain homotopy maps.
This general upper estimate is used in our construction of
spectral invariants in [Oh7].

\bigskip
\n{\it 3.2. Adiabatic homotopy and adiabatic chain map}
\smallskip

For our purpose of using the Floer theory in the study of Hofer's
geometry, we also need to consider a {\it family version } of the
Floer homology to keep track of the behavior of the action
spectrum over one parameter family of Hamiltonians as in \S2.

Let $\phi \in {\cal H}am (M,\omega)$ and $\FF = \{F^s\}_{s \in
[0,1]}$ be a path in $\HH(\phi)$. We normalize $F^s$ so that (2.5)
(and so Proposition 2.3) holds. With this normalization, {\it if $\text{\rm
Spec }(G) \subset \R$ were isomorphic to $\Gamma_\omega\Z$ or
$\{0\}$} like the case where $\pi_2(M) = 0$ or more generally
where $(M,\omega)$ is integral, the ``adiabatic'' homotopy
$$
h^{adb}_{\FF}: \widetilde {CF}(J,G) \to \widetilde {CF}(J,F)
$$
as defined in [MO1,2] will induce an isomorphism
$$
h^{adb}_{\FF}: \widetilde {HF}^{(-\infty,a]}(J,G) \to \widetilde
{HF}^{(-\infty,a]}(J,F)
$$
for any $a \in \R$. Since we will use this adiabatic homotopy in
an essential way later, we carefully explain  how it is
constructed following the exposition from [MO1,2].

Suppose that there is a `gap' in the spectrum $\text{Spec }(G) =
\text{Spec }(F^s)$, i.e, that there is a positive number $\e >0$
such that
$$
|\lambda - \mu| \geq \e
$$
for all $\lambda \neq \mu \in \text{Spec }(G)$.

Since $s\mapsto F^s$ is a smooth path, there exists some $\delta
>0$ such that
$$
\| F^u-F^s\|_{C^0} <{\e\over 3}\tag 3.7
$$
for all $u,s\in [0,1]$ with $|u-s|<\delta$.  We consider the
partition
$$
I: 0=t_0<t_1<\cdots <t_N=1
$$
so that
$$
|t_j-t_{j+1}|<\delta\;\;\hbox{for all $j$}.
$$
By Proposition 3.2, the chain map
$$
h_{us}^{lin}:\widetilde{CF}(F^u)\rightarrow \widetilde {CF}(F^s)
$$
over the linear path
$$
\LL: r \mapsto (1-r)F^u + r F^s; \, r \in [0,1]
$$
restricts to
$$
h_{us}^{lin}:\widetilde {CF}^{(-\infty ,\lambda]}(F^u)\rightarrow
\widetilde {CF}^{(-\infty ,\lambda +{\e\over 3}]}(F^s)
$$
for any $u,s\in [0,1]$ with $|u-s|<\delta$.  Similarly, we have
$$
h_{su}^{lin}:\widetilde {CF}^{(-\infty ,\lambda']}(F^s)\rightarrow
\widetilde {CF}^{(-\infty , \lambda'+{\e\over 3}]}(F^u)
$$
for any $\lambda^\prime \in \R$. Combining these two, we have the
composition
$$
h_{su}^{lin}\circ h_{us}^{lin}:\widetilde {CF}^{(-\infty
,\lambda]}(F^u)\rightarrow \widetilde {CF}^{(-\infty ,\lambda
+{2\e\over 3}]}(F^u).
$$
By the condition (3.5) and the gap condition, all of these three
maps in fact restrict to the same levels and induces homomorphisms
$$
\aligned
 &h_{us}^{lin}:\widetilde{HF}^{(-\infty
,\lambda]}(F^u)\rightarrow \widetilde{HF}^{(-\infty
,\lambda]}(F^s)\\
&h_{su}^{lin}:\widetilde{HF}^{(-\infty ,\lambda]}(F^s)\rightarrow
\widetilde{HF}^{(-\infty ,\lambda]}(F^u)
\endaligned
\tag 3.8
$$
and
$$
h_{su}^{lin}\circ h_{us}^{lin}:\widetilde{HF}^{(-\infty
,\lambda]}(F^u)\rightarrow \widetilde{HF}^{(-\infty
,\lambda]}(F^u),
$$
provided $\lambda$ is chosen sufficiently close to $\hbox{\rm Spec
}(G)$. However, if we choose $\delta$ sufficiently small, we can
also prove the identity
$$
h_{su}^{lin}\circ h_{us}^{lin}=h_{uu} (=\hbox{id })\; \text {\rm
on}\;\; \widetilde{HF}^{(-\infty ,\lambda)}(F^u)
$$
 which implies that (3.8) is an
isomorphism for all $u,s$ with $|u-s|<\delta$. By repeating the
above to $(u,s)=(t_j, t_{j+1})$ for $j=0,\ldots ,N-1$, we conclude
that the composition
$$
h_{t_jt_{j-1}}^{lin}\circ h_{t_{j-1}t_{j-2}}^{lin}\circ
\cdots\cdot\circ h_{t_1t_0}^{lin}:\widetilde {CF}(G)\rightarrow
\widetilde {CF}(F^{t_j})\tag 3.9
$$
restricts to
$$
h_{t_jt_{j-1}}^{lin}\circ\cdots\circ h_{t_1t_0}^{lin}:\widetilde
{CF}^{(-\infty ,\lambda]}(G)\rightarrow \widetilde {CF}^{(-\infty
,\lambda]}(F^{t_j})
$$
for all $1 \leq j \leq N$, and so induces the composition
$$
h_{t_jt_{j-1}}^{lin}\circ\cdots \circ
h_{t_1t_0}^{lin}:\widetilde{HF}^{(-\infty ,\lambda]}(G)\rightarrow
\widetilde{HF}^{(-\infty ,\lambda]}(F^{t_j})\tag 3.10
$$
which becomes an isomorphism. In particular, we have the
isomorphism
$$
h_{t_Nt_{N-1}}^{lin}\circ\cdots \circ
h_{t_1t_0}^{lin}:\widetilde{HF}^{(-\infty ,\lambda]}(J,
G)\rightarrow \widetilde{HF}^{(-\infty ,\lambda]}(J, F)\tag 3.11
$$

\medskip
\n{\bf Definition 3.4.} ~ Let $ I: 0 < \e_1 < \e_2 < \cdots < \e_N
=1 $ be a partition. We define its {\it mesh}, denoted as
$\Delta_I$, by
$$
\Delta_I : = \max_{j} |t_{j+1} -t_j|.
$$
We call the associated piecewise continuous linear path $\LL_1\#
\LL_2\# \cdots \# \LL_{N-1}$ and the chain map (3.9) the {\it
adiabatic homotopy}, denoted as $\FF_I$, and the {\it adiabatic
chain map over the path $\FF$}. We denote
$$
h_\FF^I = h_{t_jt_{j-1}}^{lin}\circ\cdots \circ
h_{t_1t_0}^{lin}:\widetilde{CF}(G)\rightarrow \widetilde{CF}(F).
\tag 3.12
$$
We define the {\it mesh } $\Delta(\FF_I)$ of the adiabatic
homotopy $\FF_I$ along the path $\FF$ to be
$$
\Delta(\FF_I): = \max_{j=0, \cdots, N-1}\Big\{ \int_0^1
-\min(F^{t_{j}} - F^{t_{j+1}})\, dt, \,\int_0^1 \max(F^{t_{j}} -
F^{t_{j+1}})\, dt \, \Big\} \tag 3.13
$$
We simply denote by $\FF^{adb}$, $h_{\FF}^{adb}$ when we do not
specify the partition $I$. Note that the mesh of the adiabatic
homotopy can be made arbitrarily small by making $\Delta_I$ small.

This adiabatic construction of homotopy in the chain level will be
used in a crucial way to study the global case of length
minimizing property of geodesics, where the action spectrum is not
necessarily fixed and does not have a `gap' in general.

\medskip

\head\bf \S4. $C^2$-small Hamiltonians and local Floer complex
\endhead

\n {\it 4.1. Local Floer complex}
\smallskip

In this section, we consider $C^2$-small Hamiltonians $F$. We
consider the subset $\Omega_{N_\Delta}(M)$ of loops $\gamma$ with
$(\gamma(0), \gamma(t)) \in M\times M$ is contained in a fixed
Darboux neighborhood $N_\Delta$ of the diagonal $\Delta \subset
M\times M$ for all $t \in [0,1]$. In particular, any periodic
orbit $z$ of $X_H$ contained in $\Omega_{N_\Delta}(M)$ has a
canonical isotopy class of contraction $w_z$. We will always use
this convention $w_z$ whenever there is a canonical contraction of
$z$ like in this case of small loops. This provides a canonical
embedding of $\Omega_{N_\Delta}(M)\subset \widetilde\Omega_0(M)$
defined by
$$
z \to [z,w_z].
$$

We denote by $\HH_\delta$ the set
$$
\HH_\delta =\{ F : [0,1] \times M \to \R ~|~
\|F\|_{C^2} \leq \delta \, \text { and} \, \int_M F_t = 0 \, 
\text{ for all}\, t\}.
$$

\medskip

Imitating the construction from [Fl2] and [Oh1], we define
\medskip

\n{\bf Definition 4.1.} ~ For any $(J,F) \in \JJ_\omega(M) \times
\HH_\delta$ and for the given Darboux neighborhood $N_\Delta$ of
the diagonal $\Delta \subset M\times M$ such that
$$
\phi_F^t(\Delta) \subset \hbox{\rm Int }N_\Delta,
$$
we define
$$
\MM(J,F: N_\Delta) = \{u \in \MM(J,F) ~|~ (u(\tau)(0), u(\tau)(t))
\in \hbox{\rm Int } N_\Delta  \, \hbox{\rm for all }\, \tau \}.
\tag 4.1
$$
Consider the evaluation map
$$
ev:\MM(J,F:N_\Delta) \to \Omega_{N_\Delta}; \quad ev(u) = u(0).
$$
For each open subset $\UU \subset M\times M$ with
$\Delta \subset \UU \subset M\times M $,
we define the {\it local Floer complex} in $\Omega_{\UU}$ by
$$
\SS(J,F: \UU) : = ev(\MM(J,F: \UU)) \subset \Omega_{\UU}. \tag 4.2
$$
We say $\SS(J,F:\UU)$ is {\it isolated} in $\UU$ if its closure is
contained in the interior of $\overline{\Omega_{\UU}}$.
\medskip

The following can be proved by the same method as that of [Fl2]
(See Proposition 3.2 [Oh1]), to which we refer readers for its
proof.

\proclaim\nofrills{Proposition 4.2. }~ If $\SS(J,F:\UU)$ is
isolated in $\UU$, then for all $(J^\prime, F^\prime)$
$C^\infty$-close enough to $(J,F)$ in the $C^\infty$-topology,
$\SS(J^\prime, F^\prime:\UU)$ is also isolated in $\UU$.
\endproclaim

Using this proposition, we can define the {\it local Floer
homology}, denoted by $HF(J,F:\UU)$. Furthermore, the restriction
of the action functional to the image of the embedding
$\Omega_{N_\Delta}(M) \subset \widetilde\Omega_0(M)$ provides a
filtration on the local Floer complex. Proof of the following
proposition is standard combining existing methods in the Floer
theory (see [\S 3, Oh1]).

\proclaim\nofrills{Proposition 4.3. }~ Let $\UU$ be as above and
$F, \, F^\prime \in \HH_\delta$. Assume that $\delta > 0$ so small
that (4.2) holds for $F, \, F^\prime$. Then there exists a
canonical isomorphism, we have
$$
h_{(J:\UU)}: HF(J,F: \UU)  \to HF(J,F^\prime:\UU) \tag 4.3
$$
whose matrix elements are given by the number of solutions of
(4.4) below whose images are contained in $\UU$:
$$
\cases \frac{\part u}{\part \tau} +
J \Big(\frac{\part u}{\part t} - X_{H^{\rho(\tau)}}(u)\Big) = 0 \\
u(-\infty) = x \in CF(F:\UU), \quad u(\infty) = y \in
CF(F^\prime:\UU)\\
w_x\# u \sim w_y.
\endcases
\tag 4.4
$$
\endproclaim

Following [Oh1], we call {\it thin} trajectories the solutions of
the Cauchy-Riemann equations defining the boundary map or the
chain map whose images are contained in $\UU$.
\medskip

\n {\it 4.2. $\text{Fix }\phi^1_G$ versus $\Delta \cap \text{graph
}\phi^1_G$: comparison of two Floer homology}\smallskip

The main goal of this sub-section is to prove that when $G$ is
$C^2$-small quasi-autonomous Hamiltonian, the minimum point $x^-$,
which corresponds to a (local) maximum point of $\AA_G$ in the
local Floer complex, is homologically essential in the local Floer
complex. There does not seem to be a direct way of proving this in
the context of Floer theory of Hamiltonian {\it diffeomorphisms}.
We will need to use the intersection theoretic version of the
Floer theory of Lagrangian submanifolds between $\Delta$ and
$\text{graph }\phi_G^1$ in the product $(M, -\omega) \times
(M,\omega)$. This kind of comparison argument has been around
among the experts in the Floer theory but never been rigorously
carried out before. As we will see below, contrary to the
conventional wisdom in the literature, {\it this comparison does
not work in the chain level but works only in the homology level}.

We now compare the local Floer homology $HF(J,G: \UU)$ of
$C^2$-small Hamiltonian $G$ and two versions of its intersection
counterparts, one  $HF_{-J\oplus J, 0}(\Delta, \text{graph
}\phi^1_G:\UU)$ and the other $HF_{-J\oplus (\phi_G)^*J, 0\oplus
G}(\Delta, \Delta:\UU)$. We will be especially keen to keep track
of filtration changes.

First we note that the two Floer complexes $\MM_{-J\oplus J,
0}(\Delta, \text{graph }\phi^1_G:\UU)$ and $\MM_{-J\oplus
(\phi_G)^*J, 0\oplus G}(\Delta, \Delta:\UU)$ are canonically
isomorphic by the assignment
$$
(\gamma(t),\gamma(t)) \mapsto
(\gamma(t),(\phi^t_G)^{-1}(\gamma)(t)).
$$
and so the two Lagrangian intersection Floer homology are
canonically isomorphic: Here the above two moduli spaces are the
solutions sets of the following Cauchy-Riemann equations
$$
\cases \frac{\part U}{\part \tau} +
(-J\oplus J) \frac{\part U}{\part t}  = 0 \\
U(\tau, 0) \in \Delta, \, U(\tau,1) \in \text{graph }\phi^1_G
\endcases
$$
and
$$
\cases \frac{\part U}{\part \tau} +
(-J\oplus (\phi_G^1)^*J) \Big(\frac{\part U}{\part t} 
- X_{0\oplus G}(U)\Big) = 0 \\
U(\tau, 0) \in \Delta, \, U(\tau,1) \in \Delta
\endcases
$$
respectively, where $U = (u_1,u_2) : \R \times [0,1] \to M\times
M$. The relevant action functionals for these cases are given by
$$
\AA_0([\Gamma, W]) = - \int W^*(-\omega \oplus \omega) \tag 4.5
$$
on $\widetilde \Omega(\Delta, \text{graph }\phi_G^1: M \times M)$
and
$$
\AA_{0\oplus G}([\Gamma, W]) = \AA_0 (\Gamma, W)- \int_0^1
(0\oplus G)(\Gamma(t), t)\, dt \tag 4.6
$$
on $\widetilde \Omega(\Delta, \Delta: M \times M)$ where we denote
$$
\Omega(\Delta, \text{graph }\phi_G^1: M \times M) = \{ \Gamma:
[0,1] \to M\times M ~|~ \Gamma(0) \in \Delta, \, \Gamma(1) \in
\text{graph }\phi_G^1\}
$$
and similarly for $\Omega(\Delta, \text{graph }\phi_G^1: M \times
M)$. Again the `tilde' means the  covering space which can be
represented by the set of pairs $[\Gamma, W]$ in a similar way
(see [\S 2, FOOO] for the complete discussion on this set-up for
the Lagrangian intersection Floer homology theory). The relations
between the action functionals (4.5), (4.6) and (2.1) are evident
and respects the filtration.

Next we will attempt to compare $HF(J,G;\UU)$ and $HF_{-J\oplus J,
0\oplus G}(\Delta, \Delta:\UU)$. Without loss of generality, we
will concern Hamiltonians $G$ such that $G\equiv 0$ near $t = 0,\,
1$, which one can always achieve by perturbing $G$ without
changing its time-one map (See Lemma 5.2).

It turns out that {\it there is no direct way of identifying the
corresponding Floer complexes between the two}.

As an intermediate case, we consider the Hamiltonian $G^\prime: M
\times [0,1]$ defined by
$$
G^\prime(x,t) = \cases 0 \quad  & \text{for }\, 0\leq t \leq \frac1{2}\\
2 G(x,2t)  & \text{for } \, \frac1{2} \leq t \leq 1
\endcases,
$$
and the assignment
$$
(u_0,u_1) \in \MM_{-J\oplus J, 0\oplus G}(\Delta, \Delta:\UU)
\mapsto v \in \MM(J, G^\prime: \UU) \tag 4.7
$$
with $v(\tau,t) := \overline{u_0}\# u_1(2\tau,2t)$. Here the map
$\overline{u_0}\# u_1: [0,2] \to M$ is the map defined by
$$
\overline{u_0} \# u_1(\tau,t) =
\cases u_0(\tau,1-t) \quad & \text{for }\, 0\leq t \leq 1 \\
u_1(\tau,t-1) & \text{for } \, 1\leq t \leq 2
\endcases
$$
is well-defined and continuous because
$$
\align
\overline{u_0}(\tau, 1) & = u_0(\tau,0) = u_1(\tau,0)\\
u_1(\tau,1) & = u_0 (\tau,1) = \overline{u_0}(\tau,0).
\endalign
$$
Furthermore near $t = 0, \, 1$, this is smooth (and so
holomorphic) by the elliptic regularity since $G^\prime$ is smooth
(Recall that we assume that $G \equiv 0$ near $t =0,\, 1$.
Conversely, any element $v \in \MM(J,G^\prime:\UU)$ can be written
as the form of $\overline{u_0}\# u_1$ which is uniquely determined
by $v$. This proves that (4.7) is a diffeomorphism from
$\MM_{-J\oplus J, 0\oplus G}(\Delta, \Delta:\UU)$ to $\MM(J,
G^\prime: \UU)$ which induces a filtration-preserving isomorphism
between $HF_{-J\oplus J, 0\oplus G}(\Delta, \Delta:\UU)$ and
$HF(J, G^\prime: \UU)$

Finally, we need to relate $HF(J,G:\UU)$ and $HF(J, G^\prime:\UU)$.
For this we note that $G$ and $G^\prime$ can be connected by a one-parameter
family $\overline{G} = \{G^s\}_{0\leq s\leq 1}$ with
$$
G^s(x,t):=
\cases 0 \quad &\text{for }\, 0\leq t \leq \frac{s}{2} \\
\frac2{1+s}G(x,\frac2{1+s}t) & \text{for } \frac{s}{2} \leq t \leq
1.
\endcases
$$
And we have
$$
\phi^1_{G^s} = \phi^1_G\quad \text{for all $s \in [0,1]$}.
$$
Noting that there are only finite number of periodic trajectories
in $CF(G:\UU)$, the ``adiabatic argument'' explained in \S3 indeed
proves that the adiabatic homomorphism
$$
h^{adb}_{J,\overline{G}}: CF(G^\prime :\UU) \to  CF(G :\UU) \tag
4.8
$$
respects the filtration and so the induced homomorphism in its
homology
$$
h^{adb}_{\overline{G}}: HF(J,G^\prime:\UU) \to  HF(J,G :\UU)
$$
becomes a filtration-preserving isomorphism.

We note that $\SS(J,0: \UU)$ is isolated in $\UU$. Therefore if
follows from Proposition 3.1 that if $\|G\|_{C^2}$ and
$\|K\|_{C^2}$ are sufficiently small, both $\SS(J, G:\UU)$ and
$\SS(J, K:\UU)$ are also isolated in $\UU$.

\smallskip

We now apply the above discussion to the $C^2$-small
quasi-autonomous Hamiltonian $G$ to prove the following
homological essentialness of the  minimum points of $G$ in the
local Floer homology. Recall from the remark in the beginning of
\S2 that the minimum of $G$ corresponds to the maximum of the
action functional and vice versa. We refer to [Definition 13.2.F,
Po2] for a formulation of the homological essentialness of the
critical point and its consequence on the existence result
[Corollary 13.2.H, Po2].

\proclaim{Proposition 4.4 } Suppose that $\|G\|_{C^2} < \delta$ with
$\delta$ so small that $\text{graph }\phi^1_G \subset \UU$ lies in
the given Darboux neighborhood of $\Delta \subset M\times M$.
Suppose that $G$ is quasi-autonomous with the unique maximum point
$x^+$ and minimum point $x^-$. Then the critical point $x^-$ is
homologically essential in $\SS(J,G: \UU)$.
\endproclaim

\demo{Proof} In the above discussion, we have shown that
$\MM(J,G^\prime :\UU)$ is diffeomorphic to $\MM_{-J\oplus J,
0\oplus G}(\Delta, \Delta:\UU) \cong \MM_{(-J \oplus J),0}(\Delta,
\text{graph } \phi^1_G:\UU)$.  We will first show that the
intersection point $(x^-,x^-)$ is homologically essential in the
latter Floer complex, which in turn will imply the homological
essentialness of $x^-$ in $\SS(J,G^\prime:\UU)$.

Identifying $\UU$ with a neighborhood of the zero section of the
cotangent bundle $T^*\Delta$, we denote by $J_\Delta$ the
canonical almost complex structure on $T^*\Delta$ associated to
the Levi-Civita connection of a given Riemannian metric on
$\Delta$. Since the image of $\MM_{(-J \oplus J,0)}(\Delta,
\text{graph } \phi^1_G:\UU)$ is isolated in $\UU$, we may perturb
$-J\oplus J$ to $J^\prime$ in $\UU$ so that $J^\prime \equiv
J_\Delta$ near the boundary of $\overline \UU$ and also that
$$
\MM_{(-J \oplus J)}(\Delta, \text{graph } \phi^1_G:\UU)=
\MM_{J^\prime}(\Delta, \text{graph } \phi^1_G:\UU).
$$
We connect $J^\prime$ and $ J_\Delta $ by a path $\widetilde J_t$
on $\UU$ so that $\widetilde J_t \equiv  J_\Delta$ near the
boundary for all $t \in [0,1]$. Noting that $T^*\Delta$ is
pseudo-convex with respect to $J_\Delta$, the two local Floer
complexes $\SS_{J^\prime}(\Delta, \text{graph } \phi^1_G:\UU)$ and
$\SS_{J_\Delta}(\Delta, \text{graph } \phi^1_G:\UU)$ can be
connected by an isolated continuation in $\UU$. Recall from
Proposition 3.2 that this continuation preserves the filtration of
Floer homology.

On the other hand $\hbox{\rm graph } \phi^1_G$ is diffeomorphic to
$\hbox{\rm graph } dS \subset \UU \subset T^*\Delta$ for a
generating function of the Lagrangian submanifold $\hbox{\rm graph
} \phi^1_G \subset T^*\Delta$, if $G$ is $C^2$-small. Moreover
$x^-$ corresponds to $(x^-,x^-)$ which is the minimum point of the
generating function $S: \Delta \to \R$. Since $M$ (and so
$\Delta$) is assumed to be compact, $(x^-, x^-)$ is homologically
essential in the Morse homology of $-S$. On the other hand,
$\MM_{J_\Delta}(\Delta,\hbox{\rm graph} (-dS):\UU)$ is
diffeomorphic to $\MM^{Morse}(-S,g)$ (see [FOh1] for its proof),
 where $J_\Delta$ is the almost complex
structure on $T^*\Delta$ that is associated to the Levi-Civita
connection of a chosen metric $g$ on $\Delta$. Therefore
$(x^-,x^-)$ is homologically essential in $\MM_{J_\Delta
}(\Delta,\hbox{\rm graph } \phi^1_G:\UU)$. Combining all these, we
derive that the constant solution $x^-$ is homologically essential
in the local Floer complex $\MM(J,G^\prime:\UU)$.

By the uniqueness of the minimum points, under the chain
isomorphism (4.8), the image $h^{adb}_{\overline G}(x^-)$ must
involve $x^-$ in its expression and so $x^-$ is also homologically
essential in $\SS_{J_\Delta }(J, G:\UU)$. We refer readers to the
proof of this kind of result in a more difficult context in \S7.
\qed\enddemo

\head\bf \S5. Calculation
\endhead
In this section, we start with the proof of Theorem I in the
introduction.

We consider {\it the rescaled Hamiltonians
$$
\e\, G^\e = \e\, G(\cdot, \e t) \quad 0 < \e \leq 1.
$$
and choose $\e_0 > 0$ so small that it has no non-constant 
contractible periodic orbit
for all $0 < \e \leq \e_0$.  

We first prove the following simple lemma.

\proclaim\nofrills{Lemma 5.1. } ~ Let $\{G_i\}$ be a sequence of
smooth Hamiltonians such that $G_i \to G_0$ in $C^0$-topology and
$\phi_{G_i} \to \phi_{G_0}$ in $C^0$-topology. If all $G_i$ are
length minimizing over $[a,b]$, then so is $G_0$.
\endproclaim

\demo{Proof}~ Suppose the contrary that there exists $F$ such that
$F\sim G_0$, but $\|F\| < \|G_0\|$. We choose $\delta > 0$ with
$$
\|F\| < \|G_0\| - \delta.
$$
Therefore
$$
\|F\| \leq \|G_i\| - {1 \over 2} \delta \tag 5.1
$$
for sufficiently large $i$. We consider the Hamiltonian $F_i$
defined by
$$
\align
F_i & : = (G_i \# \overline{G_0})\# F \\
& = G_i - G_0(\phi^t_{G_i}) + F(\phi^t_{G_0}\circ
(\phi^t_{G_i})^{-1}) \tag 5.2
\endalign
$$
This generates the flow $\phi^t_{G_i}\circ (\phi^t_{G_0})^{-1}
\circ \phi^t_F$ and so $F_i \sim G_i$. This implies, by the
hypothesis that $G_i$ are length minimizing over $[a,b]$, we have
$$
\|G_i\| \leq \|F_i\|
$$
and so
$$
\|F\| \leq \|F_i\| - {1\over 2} \delta \tag 5.3
$$
for all sufficiently large $i$.
 However since $G_i \to G_0$, $\phi_{G_i} \to \phi_{G_0}$
by the hypotheses (and also so $\phi_{G_0}\circ (\phi_{G_i})^{-1}
\to id$) in $C^0$-topology, we have $F_i \to F$ in $C^0$-topology.
Therefore we have
$$
\lim_{i \to \infty} \|F_i\| \to \|F\|
$$
which gives rise to a contradiction to (5.3). \qed
\enddemo

Now, using the Floer homology theory, we would like to show
$$\|G\|  \leq \|F\|$$
for any $F \sim  G$ when the quasi-autonomous Hamiltonian $G$
satisfies the hypothesis that there is no non-constant
contractible periodic orbits. However we need to take care of a
problem before applying the Floer theory, that is, $G$ not being
time-periodic. The following lemma will be important in this
respect.

\proclaim\nofrills{ Lemma 5.2. }~ Let $H$ be a given Hamiltonian
$H: T^*M \times [0,1] \rightarrow \R$ and $\phi = \phi^1_H$ be its
time-one map. Then we can perturb $H$ so that the perturbed
Hamiltonian $H'$ has the properties

\roster \item  $\phi_{H'}^1 = \phi_H^1$

\item $H' \equiv 0$ near $t = 0$ and $1$ and in particular $H'$ is
time periodic

\item  Both $|\int^1_0  \max_x (H' - H )\ dt|$ and $|\int^1_0
\min_x (H' - H)\ dt|$ can be made as small as we want

\item If $H$ is quasi-autonomous, so is $H'$.

\item There is a canonical one-one correspondence between
$\text{\rm Per}(H)$ and $\text{\rm Per}(H')$ with their actions
fixed.

\endroster
Furthermore, this modification is canonical with the ``smallness''
in (3) can be chosen uniformly over $H$ depending only on the
$C^0$-norm of $H$.
\endproclaim

\demo{Proof} We first reparameterize $\phi^t_H$ in the following
way: We choose a smooth function $\zeta : [0,1] \rightarrow [0,1]$
such that
$$
\zeta (t) = {\cases
             0  & \text {for }\, 0 \leq t \leq {\epsilon \over 2}\\
             1  &  \text {for }\, 1 - {\epsilon\over 2} \leq t \leq 1
             \endcases}
$$
and
$$
\zeta'(t) \geq 0 \quad \hbox {for all} \quad t \in [0,1],
$$
and consider the isotopy
$$
\psi^t : = \phi_H^{\zeta (t)} .
$$
It is easy to check that the Hamiltonian generating the isotopy
$\{ \psi^t \}_{0 \leq t \leq 1} $ is $H' = \{ H'_t \}_{0 \leq t
\leq 1}$ with $H'_t = \zeta '(t) H_{\zeta (t)}$.  By definition,
it follows that $H'$ satisfies (1) and (2). For (3), we compute
$$
\align
    \int^1_0 \max_x & (H' - H) dt = \int^1_0 \max_x (\zeta ' (t)
H_{\zeta(t)} - H_t ) dt \\
       &  \le \int^1_0
\max_x \Big( \zeta ' (t) (H_{\zeta (t)} - H_t )\Big) dt +
                  \int^1_0 \max_x \Big(( \zeta ' (t) - 1) H_t\Big)dt
\endalign
$$
For the first term,
$$
\align \int^1_0 \max_x & \Big(\zeta '(t) (H_{\zeta (t)} -
H_t)\Big)dt =\int^1_0 \zeta '(t) \max_x (H_{\zeta (t)} - H_t)dt\\
& \leq \int^1_0 \zeta '(t) \max_{x,t} | H_{\zeta(t)} - H_t | dt =
\max_{x,t} |H_{\zeta (t)} (x) - H_t (x) |
\endalign
$$
which can be made arbitrarily small by choosing $\zeta$ so that
$\| \zeta  - t \|_{C^0}$ become sufficiently small. For the second
term,
$$
\align
   \int^1_0 \max_x \Big(( \zeta ' (t) - 1) H_t\Big) dt
               & \leq \int^1_0 | \zeta '(t) -1 | dt \cdot \max_{x,t}
H(x,t) \\
               & \leq \| H\|_{C^0} \int^1_0 | \zeta ' (t) - 1
|dt. \endalign
$$
Again by appropriately choosing $\zeta$, we can make
$$
\int^1_0 | \zeta '(t) - 1 | dt
$$
as small as we want.  Combining these two, we have verified
$|\int^1_0 \max_x(H'-H)\, dt|$ can be made as small as we want.
Similar consideration applies to $|\int^1_0 \min_x(H' - H) \, dt|$
and hence we have finished the proof of (3). The property (4) and
naturality of this modification are evident from the construction.
(5) follows from simple comparison of corresponding actions of
periodic orbits. \qed\enddemo

We will always perform this canonical modification in the rest of
the paper whenever we would like to consider the Cauchy-Riemann
equation associated to the Hamiltonian $H$, when $H$ is not a
one-periodic Hamiltonian.

Let $F$ be an arbitrary Hamiltonian with $F \sim G$.
We want to prove $\|G\| \leq \|F\|$. Applying Lemma
5.2 to $G$ and $F$, we may assume that $G$ and $F$ are
time one periodic, allowing small errors and then getting rid of
them by taking the limit. We will postpone the proof of the
following crucial existence result to the next sections.

From now on, we will always denote by $w_y$ the constant disc $y$
for each given constant periodic orbit $y$.

\proclaim{Proposition 5.3 }~ Suppose $F \sim G$ for
sufficiently small $\e_0$ as before. Let $k$ be a Morse function
on $M$ and consider the linear homotopy
$$
L^s = (1-s) \e k + sF. \tag 5.4
$$
Then there exists $\e_1$ such that for any $0 < \e \leq \e_1$, the
continuation equation
$$
\cases \frac{\part u}{\part \tau} + J^{\rho_1(\tau)} 
(\frac{\part u}{\part t} - X_{L^{\rho_2(\tau)}}(u)) = 0 \\
u(-\infty) = y,  \quad u(\infty) = z\\
w_{y} \#u \sim w
\endcases
\tag 5.5
$$
has a solution for some $[y,w_y] \in \text{\rm Crit } (\AA_{\e
k})$ and for some $[z,w] \in \text{\rm Crit }\AA_F$ with
$$
\AA_F([z, w]) \geq \AA_F([z^-,w_{z^-}]) (= \int_0^1 - \min G
\, dt ) \tag 5.6
$$
where we recall $w_{z^-} = \widetilde h \cdot z^-$.
\endproclaim

Assuming this proposition for the moment, we proceed with the
proof of Theorem I. The following  calculation is a slight
modification used by Polterovich [Po2] in our context which will
lead to the proof of Theorem I once we prove Proposition
5.3.

We  compute
$$
\AA_F([z,w]) - \AA_{\e\, k}([y^-,w_{y^-}]) =\int^\infty_{-\infty}
\frac{d}{d\tau}\Big\{ \AA_{L^{\rho_2(\tau)}}(u(\tau), w^-\#
u(\tau))\Big\} d\tau.
$$
We have
$$
\align \frac{d}{d\tau}\Big\{ \AA_{L^{\rho_2(\tau)}}(u(\tau), & w^+\#
u(\tau))\Big\}
 = d\AA_{L^{\rho_2(\tau)}}(\frac{\part u}{d\tau})
- \rho_2^\prime(\tau) \int_0^1 (F-\e\, k)(u(\tau))\, dt \\
& =\int^1_0\omega\Big(\frac{\part u}{\part t} -
X_{L^{\rho_2(\tau)}}(u), \frac{\part u}{\part \tau}\Big)
- \rho_2^\prime(\tau) \int_0^1 (F-\e\, k)(u(\tau))\, dt \\
& = - \int^1_0\Big |\frac{\part u}{\part t}-
X_{L^{\rho_2(\tau)}}(u)\Big |^2_J
- \rho_2^\prime(\tau) \int_0^1 (F-\e\, k)(u(\tau))\, dt \\
& \leq -  \rho_2^\prime(\tau) \int^1_0 \min(F-\e k)
 \leq - \rho_2^\prime(\tau)\Big (\int^1_0 \min F 
+ \int^1_0\min \e\, k \Big)
\endalign
$$
Therefore by integrating this over $\tau$ from $-\infty$ to $\infty$, we have
$$
\AA_F([z,w]) - \AA_{\e\, k}([y^-,w_{y^-}]) \leq \int^1_0 - \min F
+ \|\e\, k\|.
$$
On the other hand, we derive
$$
\AA_F([z,w])\geq \AA_F([z^-,w_{z^-}]) = \AA_{G}([x^-,w_{x^-}]) 
= \int_0^1 - \min G\, dt
$$
from the normalization condition (2.5), (5.6) and from the fact
that $x^-$ is the fixed minimum point over $t$.
Therefore we have
$$
\int_0^1 - \min G\, dt \leq \int^1_0 -\min F + \|\e\, k\| +
\AA_{\e\, k} ([y^-,w_{y^-}]) \leq \int^1_0 -\min F + 2\e\ \|k\|
$$
By letting $\e \to 0$, we have proven
$$
\int_0^1 \min  (G) \geq  \int^1_0 \min F \tag 5.7
$$
By considering $\overline F :=- F(\phi_F^t(x),t)$ and $\overline
G$ which generate $\phi_F^{-1}$ and $\phi_G^{-1}$ respectively, we
also prove
$$
\int_0^1 \min (-G) \geq  \int^1_0 \min(- F)
$$
which is equivalent to
$$
\int_0 ^1 \max (G) \leq  \int^1_0 \max F \tag 5.8
$$
Combining (5.7) and (5.8), we have proved
$$
\|G \| \leq \|F\|.
$$
This will finish the proof of Theorem I up to the proof of
Proposition 5.3. \qed

\head \bf \S6. Handle sliding lemma
\endhead

In this section, we study an important ingredient in our proof, 
the Floer theoretic version of the `handle sliding' lemma.

Let $H$ be any time periodic Hamiltonian and consider the Cauchy
Riemann equation
$$
\frac{\part u}{\part \tau} + J \Big (\frac{\part u}{\part t} -
X_{H(u)}(u)\Big ) = 0 \tag 6.1
$$
for generic $J$. We call a solution $u$ {\it trivial} if it is
$\tau$-independent, i.e., stationary. We define
$$
\align A_{(J,H)}  := \inf\Big\{ \int \Big| {\part u\over
\part \tau}\Big |^2  ~|~ u \,& \text{ satisfies (6.1) for some }
\varepsilon \in [0,1]\\ & \text{ and is not trivial } \}. \tag
6.2\endalign
$$
The positivity of $A_{(J,H)}$ is an easy consequence of Gromov
compactness type theorem, whose proof we omit.
\medskip

We will need a family version of $A_{(J,H)}$. When there does not
occur bifurcation of periodic orbits, one can define this to be
$$
A_{(j,\HH)} = \inf_{ 0 \leq s \leq 1}A_{(J^s, H^s)}. \tag 6.3
$$
However when there does  occur bifurcation of periodic orbits,
$A_{(j,\HH)}$ could be zero, which forces us to look at another
positive constant the definition of which should be given more
subtly to make it suitable for our purpose. In introducing this
constant, we exploit the fact that in the definition of Floer's
chain homotopy map, only index zero solutions of Floer's
continuity equation (3.5) or (6.9) below enter.

We first recall that for a generic one parameter family
$\{H(s)\}_{0\leq s\leq 1}$, there are only finite number of points
$ \SS ing = \{s_1, s_2, \cdots, s_{k_1}\} \subset [0,1]$ where
there occur either birth-death or death-birth type of bifurcation
of periodic orbits (see [Lee] for a detailed proof of this). 
Furthermore at each such $s_j$, there is
exactly one bifurcation orbit $z_j$ of $\dot x = X_{H(s_j)}(x)$
for which there is a continuous family of the pair $z^+(s), \,
z^-(s)$ of periodic orbits of $\dot x = X_{H(s)}(x)$ for $|\eta -
\eta_j| < \delta$, $\delta$ sufficiently small such that
\smallskip

(1) $z^\pm(s) \to z_j$ as $s \to s_j$,\par

(2) the Conley-Zehnder indices satisfy
$$
\mu([z^+,w^+]) = \mu([z^-,w^-]) + 1 \tag 6.4
$$
where $w^+\sim w^- \# u$ for $u$ a canonical `short' cylinder
between $z^+$ and $z^-$. This latter condition makes sense because
$z^+$ and $z^-$ are close when $\delta$ is sufficiently small.
\medskip

We now prove the following important lemma

\proclaim{Lemma 6.1 }~ Let $\{H(s)\}$ be a generic one parameter
family as above. For each $s \in [0,1] \backslash \SS ing$, we
define
$$
A^0_{(J^s,H^s)} = \{ \int \Big|{\part u \over \part \tau} \Big|^2
~|~ u \, \text{satisfies (6.1), not trivial and }\, \text{Index }
u = 0 \}
$$
and
$$
A^{reg,0}_{(j,\HH)} = \inf_{s \in [0,1] \backslash \SS ing}
A^0_{(J^s, H^s)}. \tag 6.5
$$
Then $A^{reg,0}_{(j,\HH)}$ is strictly positive.
\endproclaim
\demo{Proof} Suppose the contrary that $A^{reg,0}_{(j,\HH)} = 0$,
i.e., that there exists a sequence $r_k \in [0,1] \backslash \SS
ing$ with $r_k \to r_\infty \in (0,1)$ and $u_j$ solutions of
(6.1) for $(J^{r_k},H^{r_k})$ such that
$$
\int \Big|{\part u_j \over \part \tau}\Big|^2 \to 0, \quad
\text{Index } u_j = 0. \tag 6.6
$$
Then we must have, by choosing a subsequence if necessary,
$$
r_\infty \in \SS ing
$$
and a bifurcation orbit $z_\infty$ of $\dot x =
X_{H^{r_\infty}}(x)$ such that $u_j \to z_\infty$ uniformly and so
$$
u_j(\infty), \, u_j(-\infty) \to z_\infty.
$$
Since $u_j(\pm\infty)$ are solutions of $\dot x = X_{H^{r_j}}(x)$,
they must be the pair described in (1) above (6.4) and so
$$
\text{Index }(u_j) = \mu([z^+(r_j),w^+(r_j)]) - \mu([z^-(r_j),
w^-(r_j)] = 1.
$$
But this contradicts to the index condition in (6.6) which
finishes the proof. \qed\enddemo

Again for a generic choice of $\{H^s\}$, we may assume that there
are only finitely many points $t_i \in [0,1] \backslash \SS ing$
with $i = 1, \cdots, k_2$ at which (6.1) has exactly one
non-trivial solution $u_{t_i}$ that has Fredholm index 0. (See
[Fl1] for this kind of generic argument.) We denote
$$
\NN t = \{t_i\}_{i =1, \cdots, k_2} \subset [0,1] \backslash \SS
ing.
$$

Next we define
$$
A^{sing}_{(j,\HH)} = \min_k \{ A_{(J^{s_k}, H^{s_k})} ~|~ s_k \in
\SS ing\} \tag 6.7
$$
which is again positive by Gromov type compactness theorem. Now we
have the following crucial definition of a family version of the
constant $A_{(J,H)}$ suitable for our purpose.

\definition{Definition 6.2 } We define
$$
A^0_{(j,\HH)} = \min \{ A^{reg,0}_{(j,\HH)},\, A^{sing}_{(j,\HH)}
\} > 0.
$$
\enddefinition

The following proposition is an important ingredient of our proof.

\proclaim\nofrills{Proposition 6.3. (Handle sliding lemma) } ~ Let
$j = \{J^\eta\}$ be a (two parameter) family of almost complex
structures and $\{H(\eta)\}_{0\leq \eta \leq 1}$ be a generic
family of Hamiltonians. Let $A^0_{(j,\HH)} > 0$ be the constant
defined in Definition 6.3 and let $\eta_1, \eta_2 \in [0,1]$. Then
there exists a $\delta_0
> 0$ such that if $|\eta_1 - \eta_2| < \delta$, any finite energy
solution $u$ with
$$
\text{Index }u = 0 \tag 6.8
$$
of
$$
\frac{\part u}{\part \tau} + J^{\rho(\tau)} \Big (\frac{\part
u}{\part t} - X_{H^{\rho(\tau)}}(u)\Big ) = 0 \tag 6.9
$$
must either satisfy
$$
\int \Big| {\part u \over \part \tau}\Big|^2 \leq
\varepsilon(\delta) \tag 6.10
$$
or
$$
\int \Big| {\part u \over \part \tau}\Big|^2 \geq A_{(j,\HH)}^0 -
\varepsilon(\delta) \tag 6.11
$$
where for $\varepsilon(\delta) \to 0$ as $\delta \to 0$, provided
$\delta\leq \delta_0$. Here $H^s$ is the linear path $H^s = (1-s)
H(\eta_1) + sH(\eta_2)$ and $\rho$ is the standard function as
before.
\endproclaim

We call a solution $u$ of (6.9) {\it very short} if it satisfies
(6.10), and {\it long} if it  satisfies (6.11). We can phrase the
content of this proposition as ``Any short path is indeed very
short''.
\par

\demo{Proof of Proposition 6.3}  We prove this by contradiction.
Suppose the contrary that there exists some $\varepsilon > 0$,
$\eta_1$ and $\eta_j$ with $\eta_j \to \eta_1$ as $j \to \infty$,
and solutions $u_j$ that satisfy (6.8) and
$$
\frac{\part u_j}{\part \tau} + J^{\rho(\tau)} \Big(\frac{\part
u_j}{\part t} - X_{H_{\rho(\tau)}}(u_j)\Big) = 0 \tag 6.12
$$
but
$$
\varepsilon < \int \Big|{\part u_j \over \part \tau}\Big|^2 <
A_{(j,\HH)}^0 - \varepsilon \tag 6.13
$$
In particular, the right half  of (6.13) implies the uniform bound
on the energy of $u_j$. As $j \to \infty$, the equation (6.12)
converges to (6.1) with $H = H(\eta_1)$. By Gromov type
compactness theorem, we have a cusp curve
$$
u_\infty = \sum_k u_{\infty,k}
$$
in the limit of a subsequence where each $u_{\infty,k}$ is a
solution of (6.1) with $H = H(\eta_1)$. We also have
$$
\lim_j \int \Big|{\part u_j \over \part \tau}\Big|^2 = \sum_k\int
\Big|{\part u_{\infty, k} \over \part \tau}\Big|^2.
$$
On the other hand the left half of (6.13) implies that at least
one of $u_{\infty,k}$ is not trivial.

Now we consider three cases separately: the first is the one where
$\eta_1 \in \SS ing$ and the second where $\eta_1 \in \NN t$ and
the rest where $\eta_1 \in [0,1] \backslash (\SS ing\cup \NN t)$.
When $\eta_1 \in \SS ing$, we must have
$$
\lim_j \int \Big|{\part u_j \over \part \tau}\Big|^2 \geq
A_{(J,H(\eta_1))}^{sing} \geq A_{(j,\HH)}^0
$$
which gives rise to a contradiction to (6.13) if $j$ is
sufficiently large. On the other hand, if $\eta_1 \in \NN t$, the
cusp curve must contain one component $u_\infty$ that has Index 0
and is non-constant. Again the right hand side of (6.13) prevents this from
happening. Finally when $\eta_1 \in [0,1] \backslash (\SS ing \cup
\NN t)$, the index condition $\text{Index } u_j = 0$ and the
transversality condition implies that all components
$u_{\infty,k}$ must be constant which again contradicts to LHS of
(6.13) if $j$ is sufficiently large. This finishes the proof of
proposition. \qed\enddemo

An immediate corollary of this is the following estimate on the
action.

\proclaim\nofrills{Corollary 6.4.}~ Let $j$, $\HH$ and $\delta_0$
as in Proposition 6.3. Suppose $0 < \delta \leq \delta_0$. If $u$
is very short, then we have the lower estimate
$$
-\varepsilon(\delta) + \int_0^1 - \max_{x}(H(\eta_2) -H(\eta_1))
\,dt \leq \AA_{H(\eta_2)}(u(+\infty))  -
\AA_{H(\eta_1)}(u(-\infty)) \tag 6.14
$$
and so combined with the upper estimate (3.6), we have
$$ \align -\varepsilon(\delta) +
\int_0^1 - \max_{x}(H(\eta_2) -H(\eta_1)) \,dt & \leq
\AA_{H(\eta_2)}(u(+\infty))  -
\AA_{H(\eta_1)}(u(-\infty))\\
& \leq \int_0^1 - \min_{x}(H(\eta_2) -H(\eta_1)) \,dt. \tag 6.15
\endalign
$$
If $u$ is not very short and so must be long, then we have the
improved upper estimate
$$
\AA_{H(\eta_2)}(u(+\infty))  - \AA_{H(\eta_1)}(u(-\infty)) \leq
-A_{(j,\HH)} + \varepsilon + \int_0^1 - \min_{x}(H(\eta_2)
-H(\eta_1)) \,dt. \tag 6.16
$$
\endproclaim
\demo{Proof} A straightforward computation leads to the following
general identity
$$
\align
\AA_{H(\e_2)}(u(+\infty)) & - \AA_{H(\e_1)}(u(-\infty)) \\
& =  -\int \Big|{\part u \over \part \tau}\Big|^2_J -
\int_{-\infty}^{\infty} \rho'(\tau) \int_0^1 (H(\e_2)
-H(\e_1))(u(\tau)) \,dt \, d\tau.
\endalign
$$
Corollary 6.4 immediately follows from this and Proposition 6.3.
\qed\enddemo

We will apply the above handle sliding lemma and its corollary to
the adiabatic paths  in the next section.

\head{\bf \S 7. Non-pushing down lemma and existence}
\endhead

In this section, we will assume the main hypothesis. This is the only
section where we use the hypothesis. All the materials in other
sections are valid in arbitrary compact symplectic manifolds.

\medskip

\n{\bf Hypothesis.} Assume one of the following two cases:
\roster
\item either $(M,\omega)$ is weakly exact, i.e.,
 $\omega|_{\pi_2(M)}= 0 $ or
\item $H$ is autonomous on arbitrary $(M, \omega)$
\endroster

\medskip

In the beginning, we will approach both cases in the general
setting of quasi-autonomous cases on arbitrary $(M,\omega)$ and
then explain how non-existence of quantum contributions enter our
proof of the Non-pushing down lemma.

\medskip

\n{\bf Definition 7.1.}~ Let $H: M \times [0,1] \to \R$ be a
Hamiltonian which is not necessarily time-periodic and $\phi_H^t$
be its Hamiltonian flow. \par

\roster \item We call a point $p\in M$ a {\it time $T$ periodic
point} if $\phi_H^T(p)=p$. We call $t \in [0,T] \mapsto
\phi_H^t(p)$ a {\it contractible time $T$-periodic orbit} if it is
contractible. \par

\item When $H$ has a fixed critical point $p$ over $t \in
[0,T]$, we call $p$ {\it over-twisted} as a time $T$-periodic
orbit if its linearized flow $d\phi_H^t(p); \, t\in [0,T]$ on
$T_pM$ has a closed trajectory of period less than $1$.
\endroster
\medskip

The remaining section will be occupied by the proof of the
following result (Theorem I in the introduction).

\proclaim\nofrills{Theorem 7.2.}~  We assume one of the two
cases in the Hypothesis.
Suppose that the quasi-autonomous Hamiltonian $G$ satisfies
\smallskip

\n (i) $\phi^t_G$ has no non-constant contractible periodic orbit of
period less than one, \par

\n (ii) it has at least one fixed minimum and one fixed maximum which
are not over-twisted. 
\par
\smallskip

Then the Hamiltonian path $\phi_G^t, \, 0 \leq t \leq 1$ is length
minimizing in its homotopy class with fixed ends.
\endproclaim

\n{\bf Remark 7.3.} ~
(1) Note that the hypotheses (i) is slightly different from 
Theorem I.  However from our proof,
it will be clear that the proof for Theorem 7.2 is stable under
$C^2$-small perturbation of the Hamiltonian and so allow 
sufficiently $C^1$-small non-constant contractible periodic orbits. 
This will prove Theorem I. It is rather awkward to state how small
the perturbation can be. One might want to
consider Theorem I as a stability result of the case in Theorem 7.2.

(2) Considering $\e G^\e$ with $\e < 1$ but arbitrarily close to 1
and applying Lemma 5.2, we may assume stronger assumption ``
period less than equal to 1'' instead of ``period less than 1'' in
both (1) and (2) in the hypotheses in the theorem. We will assume
this stronger assumption in the proof.
\medskip

We consider the reparameterized Hamiltonians $\e \in [\e_0,1] \mapsto \e G^\e$.
The assumption (i) implies that there is no appearance of
non-constant contractible periodic orbit as $\e$ moves from $\e_0$
to 1. The only possible bifurcation is by that of critical points
of $\e\, G^\e$. This proves

\proclaim\nofrills{Lemma 7.4. }~ Suppose $G$ satisfies the above.
Then for each $0 <\e \leq 1$, there is one-one correspondence
between the set of contractible solutions and the set of points $x
\in M$ such that
$$
dG_t(x) = 0 \quad \text{for all } \, 0 < t \leq \e. \tag 7.1
$$
\endproclaim

\definition{Definition 7.5 } We call a point $x$ {\it
$[0,\e]$-critical point} of $G$ if $x$ satisfies (7.1). We denote
by
$$
\text{Crit}_0^\e (G)
$$
the set of $[0,\e]$-critical points of $G$.
\enddefinition

It follows from Lemma 7.4 that for any $\e' > \e \geq \e_0$ there
is a canonical injection
$$
i_{\e'\e}: \text {Crit}(\AA_{\e' G^{\e'}}) \to \text
{Crit}(\AA_{\e G^{\e}}) \hookrightarrow \text{Crit}(\AA_{\e_0
G^{\e_0}}) \tag 7.2
$$
and that there is a canonical one-one correspondence between the
set of $[0,\e]$-critical points of $G$ and that of critical points
of $\AA_{\e G^\e}$ which are of the type $[x,w_x]$. From this
description of $\text{Crit }\AA_{\e G^\e}$, it follows that there
does not emerge any new critical points of $\AA_{\e G^\e}$ as $\e$
moves from $\e_0$ to 1.

For any $[0,\e]$-critical point $x$ of $G$, we have
$$
\AA_{\e G^\e}([x,w_x]) = - \int_0^1 \e G(x, \e s) \, ds =
-\int_0^\e G(x,t)\, dt . \tag 7.3
$$
We denote
$$
\gamma_x(\e) : = -\int_0^\e G(x,t)\, dt
$$
and
$$
\gamma^\pm(\e) = -\int_0^\e G(x^\mp,t) \, dt.
$$

Using Lemma 5.1 and 5.2 and the conditions (i) and (ii) in the
statement of Theorem 7.2, by adding a small bump function around
$x^-$, we may assume, without loss of generality, that $x^-$ is
the unique minimum point of $G_t$ for each $t\in [0,1]$ and that
there is a `gap' between $-G(x^-,t)$ and $-G(x,t)$
$$
- G(x^-,t) + G(x,t) > \delta_1 \tag 7.4
$$
for all $t \in [0,1]$ for any $x \neq x^- \in \text{Crit
}_0^\eta(G)$. Similar statement holds for the maximum point $x^+$.
We will fix $\delta_1
> 0$ later in (7.15). This implies that for any $\eta \geq \e_0$
we have
$$
\AA_{\eta G^\eta}([x^-,w_{x^-}]) - \AA_{\eta G^\eta }([x,w_x]) =
\gamma^+(\eta) - \gamma_x(\eta) > \eta\delta_1 \geq \e_0\delta_1
\tag 7.5
$$
for any $[0,\eta]$-critical point $x \neq x^-$ of $G$ .

For the proof of Theorem 7.2, it will be enough to prove Proposition 5.3.
The rest of this section will be occupied by its proof.
We recall that we considered the linear homotopy $\LL =\{L^s\}$,
$$
L^s = (1-s) \e k + sF.
$$
and then studied the continuation equation
$$
\cases \frac{\part u}{\part \tau} + J (\frac{\part u}{\part t}
- X_{L^{\rho(\tau)}}(u)) = 0 \\
u(-\infty) = y^- \in \text{Crit}(K),  \quad u(\infty) = z \\
y^-\#u \sim w.
\endcases
\tag 7.6
$$
Using Lemma 5.2, after preliminary perturbation of $G$, we may
assume that there are only finitely many constant periodic
solutions of $\dot x = X_G(x)$.

We will construct a solution of the equation (7.6) in four steps:
First by considering the linear homotopy
$$
\KK: \e k \mapsto \e_0 G^{\e_0},
$$
we construct a cycle $\alpha \in (\widetilde {CF}(\e_0 G^{\e_0}),
\part_{J,\e_0\, G^{\e_0}})$ with its Floer homology class
$[\alpha]$ being non-zero, and  which is a linear combination of
the form
$$
\alpha = [x^-, w_{x^-}] + \sum_j a_j [x_j, w_{x_j}], \quad a_j \in
\Q \tag 7.7
$$
where $x_j$'s are the uniform critical points of $G_t$ over $t \in
[0,\e_0]$. This is an immediate consequence of homological
essentialness (Proposition 4.4) of $x^-$ in the local Floer
complex $CF(\e\, G^\e:\UU)$  and from the Hypothesis above,
which implies that there is no quantum contribution
for the Floer boundary operator 
for the $C^2$-small Hamiltonians in either case. 
(See Proposition 7.6 below).

Secondly we consider the homotopy
$$
\GG: \eta \mapsto \{\eta G^\eta\}, \quad \eta\in [\e_0,1]
$$
from $\e_0 G^{\e_0}$ to $G$. This step proves that the Novikov
cycle $\alpha_G$ of $G$ transferred from $\alpha$ via the adiabatic 
homotopy along $\GG$ satisfies the Non-pushing down lemma, i.e, cannot be
pushed down by the Cauchy-Riemann flow of $G$. The proof
heavily relies on the Hypothesis.

Thirdly we consider the homotopy
$$
\FF: s \mapsto \{F^s\}, \quad s \in [0,1]
$$
from $G$ to $F$ which is provided by the definition $G \sim F$.
Again this step proves that the Novikov cycle of $F$
transferred from $\alpha_G$ via the adiabatic homotopy
along $\FF$ cannot be pushed down by the Cauchy-Riemann flow of $F$.
However its proof do not use the Hypothesis but the fact $G\sim F$
and the arguments hold in general.

Finally, we glue the homotopies $\KK, \,\GG $ and $\FF$ and deform the
glued homotopy $\KK \#_{R_1}\GG \#_{R_2} \FF$ to the linear homotopy
$$
\LL: s \mapsto (1-s)\e_0 G^{\e_0} + sG.
$$
The arguments in this step are independent of the Hypothesis.
\medskip

In the rest of this section, we will carry out these steps.

\n{\it Step I; from $\e k \to \e_0 G^{\e_0}$}
\smallskip

To carry out the first step, it is essential to further analyze the
general structure of the boundary operator for the $C^2$-small 
Hamiltonians (not necessarily quasi-autonomous) like $\e\, G^\e$ 
of $\e$ sufficiently small. This will be carried out following the
argument used in [\S3, Oh1].

For each time independent $J_0$, we consider the quantity
$$
A= A(J_0,\omega:M) := \inf\Big\{ \int v^*\omega ~|~ v: S^2 \to M,
\overline\partial_{J_0} v = 0, v \, \text{non-constant} \Big\}.
$$
We choose $\e >0$ so small and in particular so that $\|\e\,
G^\e\| < \frac1{2} A(J_0,\omega:M)$.

We now state the following proposition, which is the analog of
[Proposition 4.1, Oh1] to which we refer its proof (see also [Oh7] for
its complete proof).

\proclaim{Proposition 7.6 } Let $\UU$ be the Darboux neighborhood
of $\Delta$ in $M\times M$ chosen as before. Then, for any given
$\alpha > 0$ and for any fixed time-independent $J_0$, there
exists a constant $\delta > 0$ such that $\|\e\, G^\e\|_{C^2} <
\delta$ and $|J - J_0|_{M \times [0,1]} < \delta$, we have
$$
\int \Big |\frac{\part u}{\part \tau}\Big |^2_J < A(J_0,\omega:M)
- \alpha. \tag 7.8
$$
In particular, such a path has trivial homotopy class and so
$$
\int \Big |\frac{\part u}{\part \tau}\Big |^2_J < \|\e\, G^\e  \|.
\tag 7.9
$$
Moreover, all the other $u \in \MM(J,G)$ which are not contained in
$\MM(J,G:\UU)$ satisfy
$$
\int \Big |\frac{\part u}{\part \tau}\Big |^2_J > A(J_0,\omega:M)
- \epsilon_1 \tag 7.10
$$
for sufficiently small $\epsilon_1 = \epsilon_1(\delta)$ which is
independent of $\alpha$.
\endproclaim

By the argument similar to [\S 8, Oh1], we deduce that for $(J,\e
G^\e)$ chosen as above, the boundary map
$$
\part = \part_{J,\e\, G^\e}: \widetilde{CF}(\e\, G^\e) \to 
\widetilde{CF}(\e\, G^\e)
$$
is decomposed into
$$
\part = \part_{0,\e\, G^\e} + \part^\prime_{\e\, G^\e}
\tag 7.11
$$
such that $\part^\prime_{\e\, G^\e}$ maps
$\widetilde{CF}^{(-\infty,\lambda]}(G) \to
\widetilde{CF}^{(-\infty,\lambda-A+\|\e\, G^\e \|]} $. Here the
part $\part_0$ is derived from the `thin' trajectories $u$ and
$\part^\prime_{\e\, G^\e}$ from `thick' trajectories (or from
quantum contributions). In this $C^2$-small case where
the only time-one periodic orbits are the constant ones, this `thin'
and `thick' decomposition coincides with that of 
homotopically trivial and nontrivial trajectories.
The essential point of imposing the Hypothesis is that
{\it under the Hypothesis}, $\part' = 0$ and so 
$$
\part = \part_0.
$$

Now for each given $\e \in (0,\e_1]$, we define the chain map
$$
h_{\e}^{loc}:(\widetilde{CF}(\e\, k), \part_{0,\e\, k}) \to
(CF(\e_0\, G^{\e_0}),
\part_{0,\e_0 G^{\e_0}})
$$
along the linear path
$$
\KK: s \mapsto (1-s)k + s\e_0 G^{\e_0}:= K^s
$$
by considering the equation
$$
\cases \frac{\part u}{\part \tau} + J (\frac{\part u}{\part t}
- X_{K^{\rho(\tau)}}(u)) = 0 \\
u(-\infty) = p^- , \, u(\infty)= p^+\\
w^{-}_\e\# u \sim w^{+}
\endcases
$$
for given $[p^-_\e,w^-_\e]\in \text{Crit }\AA_{\e k}$ and
$[p^+,w^+]\in \text{Crit }\AA_{\e_0 G^{\e_0}}$. The induced
homomorphisms
$$
h_{\KK}:\widetilde {CF}(J, \e\, k) \to \widetilde{CF}(J, \e_0\,
G^{\e_0})
$$
and its local version
$$
h_{\KK}^{loc}: {CF}(J, \e\, k; \UU) \to {CF}(J, \e_0\,
G^{\e_0};\UU)
$$
induces an isomorphism in homology with its inverse induced by
$h_{\KK^{-1}}$ and $h_{\KK^{-1}}^{loc}$ respectively.

Now we consider a Novikov cycle
$$
\beta = \sum a_{[p,w]} [p, w], \quad a_{[p,w]} \in \Q. 
\tag 7.12
$$

The following definition which be crucial for the minimax
argument we carry out later.
\medskip

\n{\bf Definition 7.7. }~  Let $\beta$ be a Novikov cycle in
$\widetilde {CF}(H)$. We define the {\it level} of the cycle
$\beta$ and denote by
$$
\lambda_H(\beta) =\max_{[p,w]} \{\AA_H([p,w]) ~|~a_{[p,w]}  \neq
0\, \text{in }\, (7.12) \}
$$
if $\beta \neq 0$, and just put $\lambda_H(0) = +\infty$ as usual.
\medskip

As in (7.7), we can choose a cycle $\alpha$ for
$\Big(\widetilde{CF}( \e_0 G^{\e_0}), \part_{(J,\e_0
G^{\e_0}})\Big)$
$$
\alpha = [x^-,w_{x^-}] + \sum_j a_j [x_j,w_{x_j}]
$$
with
$$
\AA_{\e_0 G^{\e_0}}([x_j,w_{x_j}]) < \AA_{\e_0
G^{\e_0}}([x^-,w_{x^-}])
$$
for all $j$, its Floer homology class satisfying $[\alpha] \neq
0$. By considering the local Floer complexes $CF(J,\e k;\UU)$ and
$CF(J, \e_0 G^{\e_0};\UU)$ and their continuation and using the
homological essentialness of the maximum point $x^-$ of $ - \e_0
G^{\e_0}$, we can write
$$
\alpha - h_\KK(\alpha_{\e k}) = \part_{\e_0 G^{\e_0}}(\gamma)
$$
for some $\gamma \in CF(\e_0 g^{\e_0};\UU)$ for each given $0 < \e
< \e_1$ so that $\alpha_{\e k}$ is a finite union
$$
\alpha_{\e k} = \sum_i^\ell a_{[p_i,w_{p_i}]} [p_i,w_{p_i}]
\tag 7.13
$$
where $p_i$'s are critical points of $k$.

\proclaim\nofrills{Lemma 7.8. }~  Assume the conditions in Theorem 7.2.
Let $\alpha$ be as above. Then for
any Novikov cycle $\beta$ homologous to $\alpha$, i.e., satisfying
$$
\alpha = \beta + \part_{\e_0 G^{\e_0}}\gamma \tag 7.14
$$
for some Novikov chain $\gamma \in \widetilde {CF}(\e_0
G^{\e_0})$, we have
$$
\lambda_{\e_0 G^{\e_0}} (\beta) \geq \lambda_{\e_0
G^{\e_0}}(\alpha). \tag 7.15
$$
\endproclaim
\demo{Proof} Note that under the main Hypothesis, we have
$$
\part_{\e_0 G^{\e_0}} = \part_{0,\e_0 G^{\e_0}}
$$
for sufficiently small $\e_0$. 
In other words, all the contributions to the boundary $\part_{\e_0
G^{\e_0}}$ come from `thin' trajectories. Since $x^-$ is the maximum
point of $-G(\cdot, t) $, there cannot be any such thin trajectory
landing at $[x^-, w_{x^-}]$.

Therefore $\beta$ must have contribution from $[x^-, w_{x^-}]$ by
(7.14) since $\alpha$ does have contribution from $[x^-,w_{x^-}]$.
Hence we must have (7.14) by definition of the level function
$\lambda_{\e_0 G^{\e_0}}$. This finishes proof of the lemma.
\qed\enddemo

\medskip

\n{\it Step II: from $\e_0G^{\e_0}$ to $G$}
\smallskip

In this step we consider the homotopy
$$
\GG: \eta \mapsto \eta G^\eta, \quad \eta \in [\e_0,1].
$$
We perturb this to a generic path
$$
\HH: \eta \mapsto H(\eta), \, \eta \in [\e_0,1]
$$
so that it satisfies the genericity condition as in the Handle
sliding lemma (See the paragraph above (6.4)). By the gap
condition and the non over-twisting condition in (ii) in Theorem 7.2,
we can continue the fixed extremum points $x^\pm$ to isolated
fixed extremum points of the perturbed path $\HH: \eta \mapsto
H(\eta)$ without having small periodic points bifurcated from
them.   In particular the perturbed path $\HH$ itself becomes
quasi-autonomous. Without loss of generality, we may assume that
these fixed extrema are the same points $x^\pm$.

Other contractible periodic orbits of $H(\eta)$ will be bifurcated
from the constant periodic orbits of $\eta G^\eta$. More
precisely, we have the following lemma.

\proclaim\nofrills{Lemma 7.9. }~ For any given $\varepsilon > 0$,
there exists a generic path $\HH: \eta \mapsto H(\eta),\, \eta \in
[\e_0,1]$ in the above sense such that for each $\eta \in [\eta_0,
1]$, for any contractible periodic orbit $z$ of $H(\eta)$ of
period one there exists $x \in \text{Crit}_0^\eta G$ such that
\par \roster

\item $$ \|z - x\|_{C^2} < \varepsilon \tag 7.16$$

\item there exists a canonical small cylinder $v$ (up to homotopy)
connecting $z$ and $x$ such that
$$
|\AA_{H(\eta)}([z,w_z]) - \AA_{\eta G^\eta}([x,w_x])| <
\varepsilon \tag 7.17
$$
and
$$
\AA_{H(\eta)}([x^-,w_{x^-}]) - \AA_{H(\eta)}([z,w_z]) > {1\over
2}\e_0\delta_1 \tag 7.18
$$
where $w_z \sim w_x \# v$.
\endroster
\endproclaim
The point of Remark 7.3 (1) is that the length minimizing property 
holds for the Hamiltonian path $\eta \mapsto H(\eta)$ which is perturbed 
from $G$ and this Hamiltonian satisfies the property assumed in Theorem I (i).
Indeed the proof below proves that this path is length minimizing.
Using Lemma 5.1, we then derive the length minimizing property of the
$G$ itself.

As in \S 3, we consider the partition
$$
I: \eta_0 = \e_0 < \eta_1 < \eta_2 < \cdots < \eta_N = 1
$$
and denote its mesh of $I$ by
$$
\Delta_I = \max_{j}|\eta_{j+1} - \eta_j|.
$$
We also consider the associated piecewise linear homotopy
$$
\HH_I:= \LL_1 \# \LL_2 \# \cdots \# \LL_{N}
$$
where $\LL_j$ is the linear homotopy
$$
s \mapsto (1-s) H(\eta_{j-1}) + sH(\eta_j).
$$
We call the above piecewise linear homotopy $\HH_I$ the {\it
adiabatic homotopy} associated to $\HH$ and the partition $I$. We
also denote the associated chain map
$$
h_{\GG_I}:= h_{\eta_N\eta_{N-1}}^{\GG, lin} \circ \cdots \circ
h_{\eta_1\e_0}^{\GG,lin} : \widetilde{CF}(H^{\e_0}) \to
\widetilde{CF}(H(1))
$$
the {\it adiabatic chain map} associated to $\HH$ and $I$. We will
just denote $\HH^{adb}$ and $h_\HH^{adb}$ respectively for the
adiabatic homotopy and the adiabatic chain map associated to $\HH$
when we do not specify the partition $I$.

Now we choose $I$ with $\Delta_I$ so small that
$$
\align \Delta(\HH_I), \,\, \Delta(\HH_I^{-1}), \,\, \Delta_I\cdot
\|H\|_{C^0}
& < {1 \over 6}\e_0 \delta_1 \tag 7.19\\
\int_0^1|H(\eta_{j+1}) (x^-,t)- H(\eta_j)(x^-,t)|\, dt & < {1
\over 6}\e_0\delta_1. \tag 7.20
\endalign
$$

We recall the Handle sliding lemma, Proposition 6.3, applied to
our perturbed family $\HH$. It is easy to see from definition that
we have
$$
A^0_{(j,\HH)} \geq {3\over 4} A_{(j,\GG)} \tag 7.21
$$
if $\HH$ is sufficiently $C^\infty$-close to $\GG$, where the
constants $A^0_{(j,\HH)}$, $A_{(j,\GG)}$ are defined as in (6.3)
and (6.5). Because there does not occur bifurcation of
contractible periodic orbits along the family $\GG$, a Gromov
compactness type argument proves $A_{(j,\GG)}
> 0$. We now state a version of Handle sliding lemma that we need
in our proof.

\proclaim\nofrills{Proposition 7.10.} ~ Let $\GG$ and $\HH$ be as
above and $j =\{J^\eta\}$ be a smooth periodic (two parameter)
family of compatible almost complex structures. Let $\eta < \eta'
\in [0,1]$. Then for any fixed $j$ and for any $\varepsilon
> 0$, there exists a constant $\delta
> 0$ such that if $0 \leq \eta' - \eta < \delta$, any finite energy
solution of
$$
\frac{\part u}{\part \tau} + J^{\rho(\tau)} \Big (\frac{\part
u}{\part t} - X_{H^{\rho(\tau)}}(u)\Big ) = 0 \tag 7.22
$$
must be either satisfies
$$
\int \Big| {\part u \over \part \tau}\Big|^2 \leq \varepsilon \tag
7.23
$$
or
$$
\int \Big| {\part u \over \part \tau}\Big|^2 \geq {1\over 2}
A_{(j,\GG)}. \tag 7.24
$$
Here $H^s$ is the linear path $H^s = (1-s) H(\eta_1) + sH(\eta_2)$
and $\rho$ is the standard function as before.
\endproclaim

By choosing $\delta_1$ and then $\Delta_I$ sufficiently small, we
will also make the constant $A_{(j,\GG)}$, satisfy
$$
A_{(j,\GG)} > 3 \delta_1 \tag 7.25
$$
which is possible because $A_{(j,\GG)}$ depends only on $\e_0$ and
$\GG$ but independent of $\delta_1$.

Next we consider the cycle
$$
\alpha_{H(1)}:= h_{\HH}^{adb}(\alpha) \tag 7.26
$$
and prove the following proposition, where the condition of no
quantum contribution enters. 

\proclaim\nofrills{Proposition 7.11. (Non-pushing down lemma II)} ~ 
Let $G$ and $(M,\omega)$
as in Theorem 7.2. Then the cycle $\alpha_{H(1)}$ has the
properties
\smallskip

\n (1) $\lambda_{H(1)} (\alpha_{H(1)}) = - \int_0^1 H(1)(x^-,t)\,
dt $
\par

\n (2) Non pushing-down lemma for $\alpha_{H(1)}$ holds, i.e., for
any Novikov cycle $\beta \in \widetilde{CF}(H(1))$ homologous to
$\alpha_{H(1)}$, we have
$$
\lambda_{H(1)}(\beta) \geq \lambda_{H(1)}(\alpha_{H(1)}).
$$
\endproclaim
\demo{Proof} We consider the family of cycles
$$
\alpha_j = h_{\eta_j\eta_{j+1}}^{\HH,lin}\circ\cdots \circ
h_{\eta_1\e}^{\HH,lin}(\alpha) \in \widetilde{CF}(H(\e_{j+1}))
$$
for $j = 0, \cdots, N-1$. We will prove the following properties
of the cycle $\alpha_j$ by induction on $j$:
\smallskip

\n (P1.j) $\alpha_j$ gets non-trivial contribution from
$[x^-,w_{x^-}] \in \text{Crit } \AA_{H(\eta)}$, \par

\n (P2.j) its level satisfies
$$
\lambda_{H(\eta_j)}(\alpha_j) = -\int_0^{\eta_j}
H(\eta_j)(x^-,t)\, dt,
$$

\n (P3.j) Non pushing down lemma for $\alpha_j$ holds, i.e., for
any Novikov cycle $\beta_j$ homologous to $\alpha_j$, we have
$$
\lambda_{H(\eta_j)}(\beta_j) \geq \lambda_{H(\eta_j)}(\alpha_j)
\tag 7.27
$$
\medskip

Once we prove this, Proposition 7.11 will follow by putting $j =
N-1$.

For $j= 0$, (P1), (P2) follow from the definition of $\alpha$ and
(P3) follows from Lemma 7.8. Now suppose (P1-3.j) hold for $j$ and
we will prove them (P1-3.j+1). We first prove (P1.j+1) and
(P2.j+1). We note that
$$
h_{\eta_{j+1}\eta_j}^{\HH^{-1},lin}\circ
h_{\eta_j\eta_{j+1}}^{\HH,lin}(\alpha_j)
$$
is homologous to $\alpha_j$ and so by (P3.j), we have
$$
\lambda_{H(\eta_j)}(h_{\eta_{j+1}\eta_j}^{\HH^{-1},lin}\circ
h_{\eta_j\eta_{j+1}}^{\HH,lin}(\alpha_j)) \geq
\lambda_{H(\eta_j)}(\alpha_j) = -
\int_0^{\eta_j}H(\eta_j)(x^-,t)\, dt.
$$
Therefore (7.19) and (P2.j) together with the upper estimate imply
$$
\lambda_{H(\eta_j)}(h_{\eta_{j+1}\eta_j}^{\HH^{-1},lin}\circ
h_{\eta_j\eta_{j+1}}^{\HH,lin}(\alpha_j))
-\lambda_{H(\eta_{j+1})}(h_{\eta_j\eta_{j+1}}^{\HH,lin}(\alpha_j))
\leq {1 \over 6}\e_0 \delta_1
$$
and so
$$
\lambda_{H(\eta_{j+1})}(h_{\eta_j\eta_{j+1}}^{\HH,lin}(\alpha_j))
\geq - \int_0^{\eta_j} H(\eta_j)(x^-,t)\, dt - {1\over 6}\e_0
\delta_1. \tag 7.28
$$
This together with Proposition 7.10 and by
(7.4), also implies that any trajectory starting from the cycle $\alpha_j$ that
lands at the critical point realizing the level
$\lambda_{H(\eta_{j+1})}(\alpha_{H\eta_{j+1}})$ must be very
short: for not very short path $u$ staring from $[z,w] \neq
[x^-,w_{x^-}]$ a generator of $\alpha_j$, it follows  from (7.24)
$$
\AA_{H(\eta_{j+1})}(u(\infty)) - \AA_{H(\eta_j)}(u(-\infty)) \leq
- {1\over 2} A_{(j,\GG)} \leq - {3\over 2}\delta_1
$$
and so
$$
\align
\AA_{H(\eta_{j+1})}(u(\infty)) 
& \leq
\AA_{H(\eta_j)}(u(-\infty)) - {3\over 2}\delta_1 \\
& \leq \AA_{H(\eta_j)}([x^-,w_{x^-}]) +{1 \over 2}\e_0 \delta_1
-{3\over 2}\delta_1 \\
& \leq -\int_0^{\eta_j} H(\eta_j)(x^-,t)\, dt +{1 \over 2}\e_0 \delta_1
-{3\over 2}\delta_1 \tag 7.29
\endalign
$$
Here the last inequality follows from (7.19), (7.20) and (P2.j). Therefore
it follows from (7.28) that such trajectory $u$ cannot 
land at a critical point that realizes
the level of $\alpha_{j+1}$ since
$$
{3\over 2}\delta_1 - {1 \over 2}\e_0\delta_1 > {1 \over 6}\e_0\delta_1.
$$
Because of (7.18) and the upper estimate, it follows that any
generator $[z,w]$ with $[z,w] \neq [x^-,w_{x^-}]$ cannot land at
the critical point of $\AA_{H(\eta_{j+1})}$ that realizes the
level of $\alpha_{j+1}$. This proves that the only possible path
realizing the level of $\alpha_j$ is a very short path $u$ such
that
$$
u(-\infty) = [x^-,w_{x^-}], \, u(\infty) = [x^-,w_{x^-}].
$$
This prove (P1.j+1) and (P2.j+1).

Now it remains to prove (P3.j+1). We prove this by contradiction.
Suppose that there is a Novikov cycle $\beta_{j+1} \in
\widetilde{CF}(H(\eta_{j+1}))$ homologous to $\alpha_{j+1}$i.e.,
$$
\alpha_{j+1} = \beta_{j+1} + \part \gamma_{j+1} \tag 7.30
$$
 but
$$
\lambda_{\eta_{j+1}G^{\eta_{j+1}}}(\beta_{j+1}) <
\lambda_{\eta_{j+1}G^{\eta_{j+1}}}(\alpha_{j+1}). \tag 7.31
$$
We study the two cases separately: \roster
\item where $(M,\omega)$ is weakly exact
\item where $G$ is autonomous.
\endroster

In the case where $(M,\omega)$ is weakly exact, (7.31) indeed
implies
$$
\lambda_{H(\eta_{j+1})}(\beta_{j+1}) <
\lambda_{H(\eta_{j+1})}(\alpha_{j+1}) - {1 \over 2} \e_0 \delta_1
\tag 7.32
$$
by (7.18) because action depends only on $z$ not on the choice of
$w$. Then the upper estimate and (7.19) and (7.20) imply
$$
\align \lambda_{H(\eta_j)}
(h_{\eta_{j+1}\eta_j}^{\HH^{-1},lin}(\beta_{j+1})) & \leq
\lambda_{H(\eta_{j+1})}(\beta_{j+1}) + {1\over 6}\e_0\delta_1 \\
& < \lambda_{H(\eta_{j+1})}(\alpha_{j+1}) - {1 \over 2} \e_0
\delta_1 + {1\over 6}\e_0\delta_1 \\
& = - \int_0^{\eta_{j+1}} H(\eta_{j+1})(x^-,t)\, dt - {1 \over 3}
\e_0
\delta_1 \\
& \leq - \int_0^{\eta_j} H(\eta_j)(x^-,t)\, dt -
\int_{\eta_j}^{\eta_{j+1}} H(\eta_{j+1})(x^-,t)\, dt  \\
& \quad  - \int_0^{\eta_j}
(H(\eta_{j+1})(x^-,t) - H(\eta_j)(x^-,t)) \, dt - {1 \over 3} \e_0\delta_1\\
& \leq - \int_0^{\eta_j} H(\eta_j)(x^-,t)\, dt +  {1 \over 3}
\e_0\delta_1 - {1 \over 3} \e_0\delta_1 \\
& =- \int_0^{\eta_j} H(\eta_j)(x^-,t)\, dt =
\lambda_{H(\eta_j)}(\alpha_j)
\endalign
$$
and hence
$$
\lambda_{H(\eta_j)}(h_{\eta_{j+1}\eta_j}^{\HH^{-1},lin}(\beta_{j+1}))
< \lambda_{H(\eta_j)}(\alpha_j).\tag 7.33
$$
However (7.33) is a contradiction to (P3.j) since the cycle
$h_{\eta_{j+1}\eta_j}^{\GG^{-1},lin}(\beta_{j+1})$ is homologous
to
$$
h_{\eta_{j+1}\eta_j}^{\GG^{-1},lin}(\alpha_{j+1}) =
h_{\eta_{j+1}\eta_j}^{\GG^{-1},lin}\circ
h_{\eta_j\eta_{j+1}}^{\GG,lin}(\alpha_j)
$$
which is in turn homologous to $\alpha_j$. This finishes proof of
(P3.j+1) for this case (1).

When $G$ is autonomous, we use a generic family of $\HH =
\{H(\eta)\}$ of {\it autonomous} Hamiltonians $H(\eta)$ which are
Morse except at a finite set of $\eta$'s, and of $j = \{J^\eta\}$
where each $J^\eta$ is $t$-independent. Since $x^-$ is the minimum
point of $H(\eta)$, there is no $t$-independent trajectory  of
$\AA_{H(\eta)}$ landing at $[x^-,w_{x^-}]$. Therefore any Floer
trajectory landing at $[x^-,w_{x^-}]$ must be
$t$-dependent. Let the trajectory start at $[x,w]$, $x \in
\text{Crit} H(\eta)$ with
$$
\mu([x,w]) - \mu([x^-,w_{x^-}]) = 1, \tag 7.34
$$
and denote by $\MM_{(J^\eta, H(\eta))}([x,w],[x^-,w_{x^-}])$ the
corresponding Floer moduli space of connecting trajectories. The
general index formula shows
$$
\mu([x,w]) = \mu([x,w_{x}]) + 2 c_1([w]). \tag 7.35
$$
We consider two cases separately:  the cases of $c_1([w]) = 0$ or
$c_1([w]) \neq 0$. If $c_1([w]) \neq 0$, we derive from (7.34),
(7.35) that $x \neq x^-$. This implies that any such trajectory
must come with (locally) free $S^1$-action, i.e., the moduli space
$$
\widehat{\MM}_{(J^\eta,H(\eta))}([x,w],[x^-,w_{x^-}]) =
\MM_{(J^\eta,H(\eta))}([x,w],[x^-,w_{x^-}])/\R
$$
and its stable map compactification have a locally free
$S^1$-action {\it without fixed points}. Therefore after
$S^1$-invariant perturbation $\Xi$ via considering the quotient
Kuranishi structure [FOn] on the quotient space  
$\widehat{\MM}_{(J^\eta,H(\eta))}([x,w],[x^-,w_{x^-}])/S^1$,
the corresponding perturbed moduli space
$\widehat{\MM}_{(J^\eta,H(\eta))}([x,w],[x^-,w_{x^-}]; \Xi)$
becomes empty. This is because the quotient Kuranishi structure
has virtual dimension -1 by the assumption (7.34). We refer to
[FOn] or [LT] for more explanation on this $S^1$-invariant
regularization process. Now consider the case $c_1([w]) = 0$. 
First note that (7.34) and (7.35) imply that $x \neq x^-$.
On the other hand, if $x\neq x^-$, the same argument as above shows 
the perturbed moduli space becomes empty. 

It now follows that there is no trajectory of index 1 that land at
$[x^-,w_{x^-}]$ after the $S^1$-invariant regularization. This
together with (7.31) gives rise to a contradiction to (7.30) as in
Lemma 7.8 and finishes the proof of (P3.j+1) for the second case
(2). Hence the proof of Proposition 7.11. \qed\enddemo

\definition{Remark 7.12 } 
(1) We would like to note that a (Morse) gradient trajectory
of the Morse function $H(\eta)$
is not necessarily regular as a Floer gradient trajectory i.e., as
a solution of the perturbed Cauchy-Riemann equation, unless the
$C^2$-norm of $H(\eta)$ is sufficiently small. The ``slowness''
condition introduced in [En], [MS] is related to this problem.

(2) A careful look of the above proofs shows
that the only obstacle to extending them to arbitrary
quasi-autonomous Hamiltonians on general symplectic manifolds is
that Non pushing-down lemma will not be available for the cycle
$$
\alpha_{H(1)} = h_\HH^{adb}(\alpha)
$$
defined in (7.26) in case quantum contribution exists for the
Floer boundary operator. This will prevent us from using the
deformation argument used in the end of \S7 to produce a solution
for the continuity equation along the linear path $\LL$. Some
simpleness condition as in [BP] enables us to prove Non-pushing
down lemma, which we will investigate further elsewhere.
\enddefinition

\n{\it Step III; from $G$ to $F$}
\smallskip

Now we consider the homotopy $\FF= \{F^s\}_{0\leq s\leq 1}$
$$
G \mapsto F.
$$
We take a partition
$$
I: 0 = s_0 < s_1 < \cdots < s_{N-1} < s_N =1
$$
and its associated adiabatic homotopy $\FF^{adb}$.

We first recall from Proposition 2.3 that
$$
\text{Spec}(F^s) = \text{Spec}(G)
$$
which is of measure zero subset $\R$. We consider the family of cycles
$$
h_{\FF^s}^{adb}(\alpha), \quad s\in [0,1]
$$
and its level function
$$
\mu(s) : = \lambda_{F^s}\Big(h_{\FF^s}^{adb}(\alpha) \Big), \quad
s\in [0,1].
$$
Here $\FF^s$ is the path $t \mapsto F^{ts}, \, t \in [0,1]$. We
will provide the proof of the following proposition in the
appendix.

\proclaim{Proposition 7.13 }~ The function $\mu$ is continuous and
so constant. In particular, the cycle
$$
\alpha_F:= h_{\FF}^{adb}(\alpha)
$$
has the level
$$
\lambda_F(\alpha_F) = \alpha_{G}(\alpha) =
\int_0^{\e_0} G(x^-,t)\, dt. \tag 7.36
$$
\endproclaim

With this proposition at out disposal, we prove 

\proclaim\nofrills{Proposition 7.14 (Non pushing-down lemma III).}~ Let
$\alpha_F$ be as above. If a Novikov cycle $\beta'$ is homologous
to $\alpha_F$ in $\widetilde{CF}( F)$, i.e., satisfies
$$
\alpha_F = \beta' + \part_F (\gamma') \tag 7.37
$$
then we must have
$$
\lambda_F(\beta') \geq  \lambda_F(\alpha_F). \tag 7.38
$$
\endproclaim
\demo{Proof} Suppose the contrary that there exists $\beta'$ and
$\gamma'$ with (7.37) and
$$
\lambda_F(\beta') < \lambda_F(\alpha_F)   \tag 7.39
$$
satisfied.  We apply the  homotopy $h_{\FF^{-1}}^{adb}$ to (7.37).
Composing this with $h_{\FF}^{adb}$, we get the identity
$$
id -h_{\FF^{-1}}^{adb}\circ h_\FF^{adb}=
\part_{G} \circ \widetilde H + \widetilde H \circ
\part_{\e_0 G^{\e_0}} \tag 7.40
$$
for the obvious Floer chain homotopy $\widetilde H: \widetilde{CF}
(\e_0 G^{\e_0}) \to \widetilde{CF}(\e_0G^{\e_0})$ in a standard
way. We apply (7.40) to the cycle $\alpha$ to get
$$
\alpha - h_{\FF^{-1}}^{adb}(\alpha_F) =
\part_{\e_0G^{\e_0}} \widetilde H (\alpha_{\e_0G^{\e_0}}) \tag 7.41
$$
from the definition of $\alpha_F$ in (7.41). Inserting (7.41) into
(7.40) and using the chain property of $h_{\FF^{-1}}^{adb}$, we get
$$
\alpha -h_{\FF^{-1}}^{adb}(\beta') =
\part_{\e_0G^{\e_0}}(\widetilde H (\alpha) + h_{\FF^{-1}}^{adb}
(\gamma')). \tag 7.42
$$
Lemma 7.8 implies that
$$
\lambda_{\e_0G^{\e_0}}(h_{\FF^{-1}}^{adb}(\beta') ) \geq
\lambda_{\e_0 g^{\e_0}}(\alpha) = \e_0 c^+ \tag 7.43
$$
On the other hand, using (7.39), (7.43) and the Handle sliding lemma, and
applying the proof of Proposition 7.13 in Appendix to $\beta'$
backwards, $F \mapsto \e_0G^{\e_0}$, we prove that the function $s
\mapsto \lambda_{FF^s}^{\FF^{-1}, adb}(\beta')$ is continuous and
so constant. In particular, we have
$$
\lambda_{\e_0G^{\e_0}}(h_{\FF^{-1}}^{adb}(\beta'))
=\lambda_F(\beta').
$$
Therefore we have proven
$$
\lambda_{\e_0G^{\e_0}}(h_{\FF^{-1}}^{adb}(\beta'))
=\lambda_F(\beta') < \lambda_F(\alpha_F)
=\lambda_{\e_0G^{\e_0}}(\alpha) \tag 7.44
$$
Now (7.43) and (7.44) give rise to a contradiction. This finishes
the proof. \qed
\enddemo

\medskip
\n{\it Step IV: from the $\KK \#_{R_1} \GG^{adb} \#_{R_2} \FF^{adb}$ 
to $\LL$}
\smallskip

Finally we consider the linear homotopy $\LL= \{L^s\}_{0\leq s\leq
1}$ from $\e k$ to $F$
$$
L^s = (1-s)\e k + sF
$$
and the associated chain map
$$
h_{\LL}: \widetilde{CF}(J^0,\e k) \to \widetilde{CF}(J^1, F)
$$
(by connecting $J^0$ and $J^1$ by a generic path $\{J^s\}$).

We connect the glued homotopy $\LL_0 = \KK \#_{R_1}\GG^{adb} \#_{R_2}
\# \FF^{adb}$ and
$\LL_1 = \LL$ by any generic homotopy (of homotopies) $\overline
\LL = \{\LL_\kappa\}_{0\leq \kappa \leq 1}$ and consider the
parameterized equation
$$
\frac{\part u}{\part \tau} + J^{\rho_1(\tau)}\Big(\frac{\part u}{\part t} 
- X_{H_\kappa^{\rho_2(\tau)}}(u)\Big) = 0
$$
for $\kappa \in [0,1]$. Again this parameterized equation induces
the identity
$$
h_\LL - h_{\KK \#_{R_1}\GG^{adb} \#_{R_2}
 \FF^{adb}} = H_{\overline \LL} \part_{(J^0,\e g)} +
\part_{(J^1,F)} H_{\overline \LL}
$$
for the corresponding chain homotopy $H_{\overline \LL}:
\widetilde{CF}(J,\e g) \to \widetilde{CF}(h^*J, F)$. Applying this
identity to $\alpha_{\e k}$ above, we have
$$
h_\LL(\alpha_{\e k}) - h_{\KK \#_{R_1}\GG \#_{R_2}
\FF^{adb}}(\alpha_{\e k})
= \part_{(J^1,F)} H_{\overline \LL} (\alpha_{\e k}).
$$
Since standard gluing theorem in the Floer theory implies
$$
h_{\KK \#_{R_1}\GG^{adb} \#_{R_2}
\FF^{adb}} = h_{\FF}^{adb}\circ h_\GG^{adb} \circ h_{\KK}
$$
for sufficiently large $R_i> 0$, we have
$$
h_{\KK \#_{R_1}\GG^{adb} \#_{R_2}
\FF^{adb}}(\alpha_{\e k}) = h_{\FF}^{adb}\circ h_\GG^{adb} \circ h_{\KK}
(\alpha_{\e k}) = h_\FF^{adb}(\alpha_{H(1)}) = \alpha_F.
$$
Obviously $h_\LL(\alpha_{\e k})$ is a Novikov cycle in $\widetilde
{CF}(J^1,F)$.  Therefore Proposition 7.14 implies that
$$
\lambda_F(h_\LL(\alpha)) \geq \lambda_F(\alpha_F).
$$
By definition of the chain map $h_\LL$ and the cycle $\alpha_{\e
k}$ in (7.13), this then implies existence of $[y, w_{y}] \in
\widetilde{CF}(J, \e \, k)$ and $[z, w]\in \widetilde{CF}(h^*J,
F)$ for which there exists a solution of the following
Cauchy-Riemann equation:
$$
\cases \frac{\part u}{\part \tau} + J^{\rho_1(\tau)} \Big(\frac{\part u}{\part t}
- X_{L^{\rho_2(\tau)}}(u)\Big) = 0 \\
u(-\infty) = y,  \quad u(\infty) = z\\
w_{y}\#u \sim w
\endcases
$$
with
$$
\AA_F([z,w]) \geq  \lambda_F(\alpha_F) = \AA_F([z^-,w_{z^-}]).
$$
This is exactly what we wanted to prove for Proposition 5.3. This
finally finishes the proof of Proposition 5.3 and hence the proof
of Theorem I.

\head{\bf \S 8. Construction of spectral invariants}
\endhead

In this section, we outline a construction of spectral invariants of
Viterbo's type [V] (more precisely, the type the author
constructed in [Oh3,5]) on arbitrary compact symplectic manifolds.
As a consequence, we also define a new invariant norm on the
Hamiltonian diffeomorphism group of arbitrary compact symplectic
manifolds. We just illustrate the main idea of the construction 
in the present paper with minimal possible sophistication 
in the presentation and refer readers [Oh7] for precise 
details of the construction.

The starting point of our construction of the
invariants will then be the fact that for any fixed generic
autonomous Hamiltonian $g$ on $M$ we have the isomorphism
$$
({CF}_*(\e g; \Lambda_\omega), \part_{\e g}) \simeq ({CM}_*(-\e
g;\Q), \part_{-\e g}^{Morse}) \otimes \Lambda_\omega \tag 8.1
$$
as a chain complex when $\e
> 0$ is sufficiently small, and the canonical isomorphism
$$
h_{\e gH}:  HF_*(\e g; \Lambda_\omega) \to HF_*(H; \Lambda_\omega)
\tag 8.2
$$
for any Hamiltonian $H$ over the Novikov ring $\Lambda_\omega$. A
natural isomorphism (8.2) is induced by the chain map
$$
h_{\e gH}: \widetilde{CF}(\e g) \to \widetilde{CF}(H) \tag 8.3
$$
over the linear path $H^s = (1-s)\e g + s H$. Here we also note
that (8.1) also induces a canonical isomorphism
$$
HF_*(\e g; \Lambda_\omega) \simeq HM_*(-\e g;\Q)\otimes
\Lambda_\omega.
$$
Here $CM_*(-\e g;\Q)$ and  $HM_*(-\e g;\Q)$ denote the Morse chain
complex and its associated homology of $-\e g$ with
$\Q$-coefficients.

By letting $\e \to 0$, we will have the corresponding limit
isomorphism
$$
h_H: H_*(M; \Q) \otimes \Lambda_\omega \to HF^*(H; \Lambda_\omega)
\tag 8.4
$$
by identifying the singular cohomology $H^*(M,\Q)$ with $HM^*(\e
g;\Q)$
by realizing its Poincar\'e dual by a Morse cycle of
$-\e g$ and then composing with the map (8.1). 
\medskip

\n{\bf Definition 8.1.} ~ Let $H$ be a given generic Hamiltonian.
For each $a\neq 0 \in H^*(M; \Q)$, we denote by $PD(a)$
its Poincar\'e dual to $a$. We consider the Floer
homology class $h_{\e gH}(PD(a)) \in HF^*(H;\Lambda_\omega)$.  We define
the {\it level} of the Floer homology class $h_{\e gH}(PD(a))$ by
$$
\rho(H;a) = \lim_{\e \to 0} \inf \{\lambda_H(\alpha) ~|~ [\alpha] 
= h_{\e gH}(PD(a)), \,
\alpha \in \widetilde{CF}(H) \, \}. \tag 8.5
$$
\medskip

Of course, a crucial task in this definition is to show that this
is well-defined, i.e, the numbers are finite, independent of the choice
of the Morse function $g$ and behave
continuously over $H$ (in $C^0$-norm). The following theorem is
the analog to [Theorem II, Oh5] which can be proved in a similar
way. However we exploit the isomorphism (8.1) in a crucial way
here.

\proclaim\nofrills{Theorem 8.2. }~ Let $H$ be a given
Hamiltonian. For each $a \neq 0 \in H^*(M;\Q)$, the number
$\rho(H;a)$ is finite and the assignment $H \mapsto \rho(H;a)$ can
be extended to $C^0(M \times [0,1])$ as a continuous function with
respect to $C^0$-topology of $H$.
\endproclaim
\demo{Proof} The proof will be the same as [Oh5] once we prove
finiteness of the value $\rho(H;a)$. 

To be more precise, we choose a Morse function $g$ on $M$ and use
the chain map (8.3). The homology class $PD(a)$ considered as a Morse
homology class of $-\e g$ defines a Floer homology class of $\e g$
which is non-zero by the fact that the Floer boundary operator
$\part _{\e g} \simeq
\part_{-\e g}^{Morse}\otimes \Lambda_\omega$. Therefore we have
$$
\inf_{[\alpha] = h_{\e gH}(PD(a))}\lambda_H(\alpha) < \infty
\tag 8.6
$$
since $h_{\e gH}(PD(a)) \neq 0$. In fact, by the same calculation as
in Proposition 3.2, we can prove
$$
\rho(H;a) \leq \int_0^1 -\min H \, dt. \tag 8.7
$$
To prove $\rho(H;a) > -\infty$, we first prove the following
lemma.

\proclaim\nofrills{Lemma 8.3. }~ We have
$$
\rho(\e g:a) \geq -\max \e g. \tag 8.8
$$
\endproclaim

\demo{Proof} Let $\gamma \in \widetilde{CF}(\e g)$ be a Novikov
cycle with $[\gamma] = PD(a)$. We write
$$
\gamma = \gamma_0 + \gamma'
$$
where $\gamma_0$ is the sum of the terms with trivial homotopy
class i.e., those of the type with $[x,w_x], \, x \in \text{Crit
}g$ and $\gamma'$ are the ones $[x,w]$ with non-trivial homotopy
class with $[w] \neq 0 \in \Gamma_\omega$. Since $\part_{\e g}$
preserves this decomposition (no quantum contribution!) and since
any $b \in H_*(M;\Q)$ can be represented by $\gamma_0$, both
$\gamma_0$ and $\gamma'$ are closed and have
$$
[\gamma_0] = b \quad \text{and }\, [\gamma'] = 0.
$$
By setting $ 0 \neq b = PD(a) = [\gamma_0]$ in the Morse homology of $\e g$, we
have $\gamma_0 \neq 0$. An easy fact from the (finite dimensional)
Morse homology theory implies
$$
\lambda_{\e g}(\gamma_0) \geq \min (-\e g) = - \max(\e g). \tag
8.9
$$
Obviously we have $\lambda_{\e g}(\gamma) \geq \lambda_{\e
g}(\gamma_0)$ which finishes the proof by (8.9). \qed\enddemo

Now we go back to the proof of Theorem 8.2. Let $\alpha \in
\widetilde {CF}(H)$ with its Floer homology class $[\alpha] = h_{\e
gH}(PD(a))$. Note that by the same calculation as in Proposition 3.2
along the linear path from $H$ to $\e g$, we have
$$
\lambda_{\e g}(h_{H \e g}^{ lin}(\alpha)) \leq \lambda_{H }(\alpha)
+ \int_0^1 - \min (\e g - H)\, dt \tag 8.10
$$
where we know $h_{H \e g}^{ lin}(\alpha) \neq 0$ because $[h_{H \e
g}^{ lin}(\alpha)] \neq 0$ since $[\alpha] \neq 0$ and $h_{H \e
g}^{lin}$ induces an isomorphism in homology. On the other hand,
let $\gamma_0$ be a representative as in Lemma 8.3 with
$[\gamma_0] = b$. Since $[\alpha] =h_{\e gH}^{lin}(\gamma_0)$, we
have
$$
[h_{H \e g}^{lin}(\alpha)] = [\gamma_0].
$$
It follows from Lemma 8.3 that
$$
\lambda_{\e g}\Big(h_{H\e g}^{lin} (\alpha)\Big) \geq -\max \e g.
\tag 8.11
$$
From (8.10), we derive
$$
\align
 \lambda_H(\alpha) & \geq  \int_0^1 \min(\e g- H)\, dt - \max
\e g \\
& = \int_0^1 -\max(H -\e g)\, dt - \max \e g.
\endalign
$$
Letting $\e \to 0$, we have proved
$$
\lambda_H(\alpha) \geq \int_0^1 - \max H \, dt \tag 8.13
$$
and then taking the infimum over $\alpha \in \widetilde{CF}(H)$
with $h_{\e g H}(PD(a)) = [\alpha]$, we derive
$$
\rho(H;a) \geq \int_0^1 -\max H \, dt \tag 8.14
$$
which in particular proves $\rho(H;a) > -\infty$.

To prove continuity of $H \mapsto \rho(H;a)$ in $C^0$-topology,
we imitate the above argument by replacing $\e g$ by another
generic Hamiltonian $K$. As in (8.10), we have
$$
\lambda_K(h_{H K}^{ lin}(\alpha)) \leq \lambda_{H }(\alpha) +
\int_0^1 - \min (K - H)\, dt. \tag 8.15
$$
We have $[h_{HK}^{lin}(\alpha)] = [h_{\e g K}(PD(a))]$ in
$HF(K;\Lambda_\omega)$ because $[\alpha] = h_{\e gH}(PD(a))$ in $HF(H;
\Lambda_\omega)$. From (8.15) and the definition of
$\rho(H;a)$, we have
$$
\rho(K;a) \leq \lambda_H(\alpha) + \int_0^1 - \min(K-H)\, dt.
$$
This proves
$$
\rho(K;a) \leq \rho(H;a) + \int_0^1 - \min(K-H)\, dt
$$
by taking the infimum of $\lambda_H(\alpha)$ over $\alpha$ with
$[\alpha] = h_{\e H}(PD(a))$. Equivalently we have
$$
\rho(K;a) - \rho(H;a) \leq \int_0^1 - \min(K-H)\, dt \tag 8.16
$$
Next we want to prove
$$
\int_0^1 -\max (K -H) \, dt \leq \rho(K;a) - \rho(H;a). \tag 8.17
$$

We apply (8.15) with $H$ and $K$ switched and $\alpha'$ with
$[\alpha'] = h_{\e gK}(PD(a))$ and get
$$
\lambda_H(h_{KH}^{ lin}(\alpha')) \leq \lambda_{K }(\alpha') +
\int_0^1 - \min (H - K)\, dt
$$
or
$$
\lambda_{K }(\alpha') - \lambda_H(h_{KH}^{ lin}(\alpha')) \geq
\int_0^1 \min (H - K)\, dt = \int_0^1 - \max (K - H).
$$
Since $[\alpha'] = h_{\e gK}(PD(a))$ and $h_{KH}^{lin}\circ h_{\e gK}$
is chain homotopic to $h_{\e gH}^{lin}$, we also have
$$
[h_{KH}^{lin}\alpha'] = h_{\e gH}(PD(a)).
$$
Therefore we derive (8.17) from this by the same argument as that
of (8.16). Combining (8.16) and (8.17), we have proved
$$
\int_0^1 -\max (K-H) \, dt \leq \rho(K;a) - \rho(H;a) \leq
\int_0^1 - \min(K-H)\, dt. \tag 8.18
$$
Now it follows from (8.18) that the function $H \mapsto
\rho(H;a)$ can be extended to $C^0(M\times [0,1])$ as a continuous
function in $C^0$-topology. This finishes the proof. \qed\enddemo

These numbers $\rho(H;a)$ will satisfy the same kind of properties
as the invariants constructed  by the author in [Oh5]. We refer to
[Oh3,5] for the statements and proofs of the properties of $\rho$
in the context of Lagrangian submanifolds on the cotangent bundle
leaving complete details to [Oh7] for the present case.

We now focus on the special cases where the corresponding class
$a$ is  the class 1 in $H^*(M;\Q)$.

\medskip

\n{\bf Definition $\&$ Theorem 8.4 [Oh7]. } ~ Let  $1$
be the identity class of $H^*(M, \Q)$. For each
given Hamiltonian $H$, we define
$$
\gamma(H) = \rho(H;1) + \rho(\overline H;1).
\tag 8.19
$$
Then we have $\gamma(H) \geq 0$, and
$$
\gamma(H) = \gamma(K) 
$$
as long as $H\sim K$. This makes $\gamma(H)$ depends only on the
equivalence class $[H]$, i.e, defines a well-defined function on
the covering space $\pi: \widetilde{\HH am}(M,\omega) \to \HH
am(M,\omega)$. Now for a given Hamiltonian diffeomorphism $\phi$,
we define
$$
\gamma(\phi) = \inf_{H \mapsto \phi}\gamma(H) = \inf_{\pi([H]) =
\phi} \gamma([H]) \tag 8.20
$$
for any  Hamiltonian diffeomorphism $H\mapsto \phi$.
\medskip

The following theorem has been proven in [Oh7] to which
we refer the readers.

\proclaim\nofrills{Theorem 8.5 [Oh7]. }~ The above function
$\gamma: \HH am(M,\omega) \to \R_+$ satisfies the following
properties: \roster
\item $\gamma(\phi) = 0$ iff $\phi = id$
\item $\gamma(\phi_1 \phi_2) \leq \gamma(\phi_1) + \gamma(\phi_2)$
\item $\gamma(\psi \circ \phi \circ \psi^{-1}) = \gamma(\phi)$ for
any symplectic diffeomorphism $\psi$.
\item $\gamma(\phi) \leq \|\phi\|$
\endroster
\endproclaim
This norm reduces to the norm Schwarz constructed in [Sc]
for the symplectically aspherical case where the norm
$\gamma$ is defined by
$$
\gamma(H) = \rho(H;1) - \rho(H;\mu) \tag 8.21
$$
where $\mu$ is the volume class in $H^*(M)$, following [V] and [Oh5]. 
The reason why the two (8.19) and (8.21) coincide is 
that in the aspherical case, we have the additional identity
$$
\rho(\overline H:1) = - \rho(H;\mu). \tag 8.22
$$
But Polterovich observed [Po3] that this latter identity fails 
in the non-exact case due to the quantum contribution. 
In fact in the non-exact case, even positivity of (8.21) seems to fail.
It turns out that our definition (8.19) in Definition 8.4
is the right one to take, which satisfies all the expected properties.
We refer readers to [Oh7] for the proof of Theorem 8.5 and for
further consequences of the spectral invariants in the study of
length minimizing property of Hofer's geodesics 
and new lower bounds for the Hofer norm of Hamiltonian 
diffeomorphisms.

\head {\bf Appendix}
\endhead

In this appendix, we prove Proposition 7.13. Since this proposition
is a general fact for arbitrary pairs  $(G,F)$ of Hamiltonians
with $G\sim F$, we gather the facts from the main part of the
paper that are needed and make this appendix self-contained.

We first recall the Handle sliding lemma. Let $H$ be any
Hamiltonian and consider the Cauchy-Riemann equation
$$
\frac{\part u}{\part \tau} + J \Big (\frac{\part u}{\part t} -
X_{H(u)}(u)\Big ) = 0. \tag A.1
$$
We call a solution $u$ {\it trivial} if it is $\tau$-independent,
i.e., stationary. We define
$$
\align A_{(J,H)}  := \inf\Big\{ \int \Big| {\part u\over
\part \tau}\Big |^2  ~|~ u \,& \text{ satisfies (A.1) and is
not trivial } \}. \tag A.2
\endalign
$$
Let $j = \{J^s\}_{0 \leq s \leq 1}$ and the family
$\HH=\{H(\eta)\}_{\eta \in [0,1]}$ be given. We define
$$
A_{(j,\HH)} = \inf_{\eta \in [0,1]} A_{(J^\eta,H(\eta))}. \tag A.3
$$
In general, this number could be zero. When it becomes positive,
we have the following  result. This is an easy version of
Proposition 6.3

\proclaim\nofrills{Proposition A.1. } ~ Let $\{H(\eta)\}_{0\leq
\eta \leq 1}$ be a smooth family of Hamiltonians and  $j =\{J^s\}$
be a smooth periodic (two parameter) family of compatible almost
complex structures. Suppose that $A_{j,\HH}$ is positive.  Let
$\eta_1, \eta_2 \in [0,1]$. Then for any fixed $j$ and for any
$\varepsilon
> 0$, there exists a constant $\delta
> 0$ such that if $|\eta_1 - \eta_2| < \delta$, any finite energy
solution of
$$
\frac{\part u}{\part \tau} + J^{\rho(\tau)} \Big (\frac{\part
u}{\part t} - X_{H^{\rho(\tau)}}(u)\Big ) = 0 \tag A.4
$$
must be either satisfies
$$
\int \Big| {\part u \over \part \tau}\Big|^2 \leq \varepsilon \tag
A.5
$$
or
$$
\int \Big| {\part u \over \part \tau}\Big|^2 \geq A_{(j,\HH)}
-\varepsilon \tag A.6
$$
Here $H^s$ is the linear path $H^s = (1-s) H(\eta_1) + sH(\eta_2)$
and $\rho$ is the standard function as before.
\endproclaim
As in Proposition 6.3, we call a solution $u$ of (A.4) very short
if it satisfies (A.4) and long if it satisfies (A.6).

\proclaim\nofrills{Corollary A.2 [Corollary 6.4, \S 6].}~ Let
$\varepsilon > 0$ be any given number. Then there exists $\delta >
0$ such that for any $\eta_1, \, \eta_2$ with $|\eta_2 -\eta_1| <
\delta$,  the following holds: if $u$ is very short, then
$$
\align -\varepsilon + \int_0^1 - \max_{x}(H(\eta_2) -H(\eta_1))
\,dt & \leq \AA_{H(\eta_2)}(u(+\infty))  -
\AA_{H(\eta_1)}(u(-\infty))\\
& \leq \int_0^1 - \min_{x}(H(\eta_2) -H(\eta_1)) \,dt. \tag A.7
\endalign
$$
If $u$ is not very short, then we have
$$
\AA_{H(\eta_2)}(u(+\infty))  - \AA_{H(\eta_1)}(u(-\infty)) \leq
-A_{(j,\HH)} + \varepsilon + \int_0^1 - \min_{x}(H(\eta_2)
-H(\eta_1)) \,dt. \tag A.8
$$
\endproclaim

We would like to apply these results to the path $\FF = \{F^s\}_{0
\leq s \leq 1}$. We first prove

\proclaim\nofrills{Lemma A.3.}~ Let $j=\{J^s\}$ be the family of
almost complex structures defined by
$$
J^s_t = (h^s_t)^*J_t.
$$
Then we have
$$
A_{(J^s, F^s)} = A_{(J,\e_0 G^{\e_0})}.
$$
In particular, we have
$$
A_{(j,\FF)} > 0. \tag A.9
$$
\endproclaim
\demo{Proof} We first note that the map
$$
x \mapsto z_x; \quad z_x(t): = h^s_t(x)
$$
and (2.4) give one-one correspondence between $\text{Per}(\e_0
G^{\e_0})$ and $\text{Per}(F^s)$ and between $\text{Crit
}\AA_{\e_0G^{\e_0}}$ and $\text{Crit }\AA_{F^s}$ respectively.
Furthermore (A.10) also provides one-one correspondence between
solution sets of the corresponding Cauchy-Riemann equations by
$$
u \mapsto u^s; \quad  u^s(\tau,t) =h^s_t(u(\tau,t)).
$$
And a straightforward calculation shows the identity
$$
\int \Big|{\part u\over \part \tau}\Big|^2_J =\int \Big|{\part
u^s\over \part \tau}\Big|^2_{J^s}
$$
which finishes the proof. \qed\enddemo

We are now ready to provide the proof of Proposition 7.13. We
choose the partition
$$
I: 0 = s_0 < s_1 < \cdots < s_N = 1
$$
so that its mesh
$$
\Delta_I(\FF) < {1\over 2}\e_0\delta_1 \tag A.10
$$
where $\Delta_I(\FF)$ is defined by
$$
\Delta_I(\FF): = \inf_j \Big\{ \int_0^1 - \min(F^{s_{j+1}}
-F^{s_j})\, dt, \int_0^1 \max(F^{s_{j+1}} -F^{s_j})\, dt \Big\}.
$$
We will prove the proposition in 3 steps: {\it finiteness, upper
estimates} and {\it lower estimates}.
\medskip
\n{\it Step 1: finiteness}
\smallskip

The finiteness of this function follows from the assumption
$[\alpha] \neq 0$ and from construction of the chain map. More
specifically, the chain map
$$
h_{\FF^s}^{adb}: \widetilde{CF}(\e_0 G^{\e_0}) \to
\widetilde{CF}(F)
$$
maps Novikov cycles to Novikov cycles and induces an isomorphism
in homology over the Novikov rings. Since we have chosen $\alpha$
so that $[\alpha] \neq 0$, we have $[h_{\FF^s}^{adb}(\alpha)] \neq
0$ and in particular $h_{\FF^s}^{adb}(\alpha) \neq 0$ for all $s$.
Hence comes the finiteness of the level of
$h_{\FF^s}^{adb}(\alpha)$, i.e, the value of $\mu(s)$.
\medskip
\n{\it Step 2: upper estimates}
\smallskip
In this step, we will prove
$$
\mu(s) - \mu(s') \leq \int_0^1 -\min_x (F^s-F^{s'})\,dt \tag A.11
$$
for $s,\, s'$ with $|s -s'| \leq \delta$ for sufficiently small
$\delta$.  This upper estimates can be proved without help of
the Handle sliding lemma.

We recall that the chain map $h_{\FF^s}^{adb}$ is defined as the
composition of chain maps $h_{s_{j+1}s_j}^{lin}:
\widetilde{CF}(F^{s_{j+1}}) \to \widetilde{CF}(F^{s_j})$ over the
linear homotopy for the partition $I$. We first consider the first
segment $[0,s_1]$. In this segment, we have
$$
h_{\FF^s}^{adb} = h_{\FF^s}^{lin}
$$
over the linear path $u \mapsto (1-u)\e_0G^{\e_0} + u F^s$.

We consider the chain map $h_{ss'}^{\FF,lin}$ which is induced by
the assignment
$$
h_{ss'}^{\FF,lin}([z^-,w^-]) = \sum_{ [z^+, w^+] \in
\text{Crit}\AA_{F^{s'}} }\#\Big(\MM_J^{\GG_1}([z^-,w^-],
[z^+,w^+])\Big)[z^+, w^+]
$$
for each $[z^-,w^-] \in \text{Crit}\AA_{F^s}$. Here $
\MM_J^{\GG_1}([z^-,w^-], [z^+,w^+])$ denotes the moduli space of
trajectories  of the Cauchy-Riemann equation
$$
{\part u\over \part \tau} + J^{\rho(\tau)}\Big({\part u\over \part
t} - X_{F^{\rho(\tau)}}(u)\Big) \tag A.12
$$
and $\#\Big(\MM_J^{\GG_1}([z^-,w^-], [z^+,w^+])\Big)$ denotes its
(rational) Euler number (see [FOn], [LT], [Ru] for the precise
meaning). In the case relevant to the chain map the moduli space
is zero-dimensional. {\it In particular, if this number is not
zero, then (A.12) has a solution.}

Assuming existence of such pair $[z^-,w^-] \in
\text{Crit}\AA_{\eta g}$ and $[z^+,w^+] \in \text{Crit}\AA_{\eta'
g}$ for the moment, we proceed with the proof. Then to every pair
$[z^-,w^-]$ and $[z^+,w^+]$ for which
$\#\Big(\MM_J^{\GG_1}([z^-,w^-], [z^+,w^+])\Big)$ is non-zero, we
have
$$
\AA_{\eta' g}(u(\infty)) - \AA_{\eta g}(u(-\infty)) \leq \int_0^1
- \min_x \,(F^{s'}- F^s)\, dt. \tag A.13
$$
Taking the maximum over $[z^-,w^-]$ among the generators of
$h_{0s}^{\FF,lin}(\alpha)$, we get
$$
\AA_{\eta' g}(u(\infty)) - \mu_1(\eta) \leq \int_0^1 - \min_x
\,(F^{s'}- F^s)\, dt. \tag A.14
$$
Since this holds for any generator $[z^+,w^+]=u(\infty)$ of
$h_{0s}^{\FF,lin}(\alpha)$, (A.14)) proves (A.11) by definition of
$\mu$.

Now it remains to prove existence of a pair $[z^-,w^-]
\in\text{Crit}\AA_{F^s}$ and $[z^+,w^+] \in
\text{Crit}\AA_{F^{s'}}$ such that
$$
\#\Big(\MM_J^{\GG_1}([z^-,w^-], [z^+,w^+])\Big) \neq 0 \tag A.15
$$
and $[z^-,w^-]$ contributes $h_{0s}^{\FF,lin}(\alpha)$ and
$[z^+,w^+]$ contributes $h_{0s}^{\FF,lin}(\alpha)$. We recall that
$$
h_{0s'}^{\FF,lin} - h_{ss'}^{\FF,lin}\circ h_{0s}^{\FF,lin} =
\part_{F^{s'}}\circ \widetilde H + \widetilde H \circ \part_{F^s}
$$
where $\widetilde H$ is defined by considering parameterized
equation induced by the homotopy (of homotopies) $\overline \LL =
\{\LL_\kappa\}_{\kappa}$ connecting the linear homotopy between
$F^0 = \e_0 G^{\e_0}$ and $F^s$ and the glued homotopy via $0
\mapsto s\mapsto s'$. However if $s$ is close to $s'$ and the
Cauchy-Riemann equation for $\LL_0$ is regular, then those
corresponding to $\LL_k$ are all regular for $0 \leq \kappa \leq
1$. Since $\widetilde H$ is defined by counting generic
non-regular solutions on $\kappa \in (0,1)$, this proves that
$\widetilde H = 0$ if $s$ is very close to $s'$. Therefore we have
$$
h_{0s'}^{\FF,lin}(\alpha) = h_{ss'}^{\FF,lin}\circ
h_{0s}^{\FF,lin}(\alpha) \tag A.16
$$
if $|s-s'| < \delta$ for sufficiently small $\delta$. By
definition of the chain map $h_{ss'}^{\FF,lin}$, there must be
such a pair of $[z^-,w^-]$ and $[z^+,w^+]$ for which (A.15) holds. 
This finishes the proof of (A.11).
\medskip

\n{\it Step 3: lower estimate}
\smallskip

This is the place where the Handle sliding lemma plays a crucial
role. We apply $h_{s's}^{\FF^{-1},lin}$ to (A.16) to get
$$
h_{s'0}^{\FF^{-1},lin}\circ h_{0s'}^{\FF,lin}(\alpha) =
h_{s'0}^{\FF^{-1},lin}\circ h_{ss'}^{\FF,lin}\circ
h_{0s}^{\FF,lin}(\alpha).
$$
Therefore $h_{s'0}^{\FF^{-1},lin}\circ h_{ss'}^{\FF,lin}\circ
h_{0s}^{\FF,lin}(\alpha)$ is homologous to $\alpha$ in
$\widetilde{CF}(\e_0 G^{\e_0})$ because
$h_{s'0}^{\FF^{-1},lin}\circ h_{0s'}^{\FF,lin}(\alpha)$ is so. By
Non pushing-down lemma, Proposition 7.14, we have
$$
\lambda_{\e_0 G^{\e_0}}(h_{s'0}^{\FF^{-1},lin}\circ
h_{ss'}^{\FF,lin}\circ h_{0s}^{\FF,lin}(\alpha)) \geq
\lambda_{\e_0 G^{\e_0}}(\alpha) = \e_0 c^+.
$$
This gives rise to
$$
\align
\lambda_{F^{s'}}(h_{ss'}^{\FF,lin}\circ h_{0s}^{\FF,lin}(\alpha))
& \geq \lambda_{\e_0 G^{\e_0}}(\alpha) - { 1\over 2} \e_0 \delta_1 \\
& \geq \lambda_{F^s}(h_{0s}^{\FF,lin}(\alpha)) - 
\e_0\delta_1.
\endalign
$$
Now we choose $\delta_1$ so small in (7.4) that we have
$$
\e_0 \delta_1 < \min\Big\{{1\over 2}A_{(j,\FF)},{1\over
2}A_{(j^{-1},\FF^{-1})}\Big\}.
$$
Then the trajectory constructed in Step 2 that satisfies (A.15)
must be very short.  On the other hand for very short
trajectories, the lower estimate (A.7) holds.

Combining Step 1-3, we have proved that the function $\mu$ is
continuous and so must be constant on $[0,s_1]$. Then this also
implies Non pushing-down lemma for $F^{s_1}$ from which we can
repeat the above argument to the segment $[s_1,s_2]$. We repeat
this to all $j = 3, \cdots, N-1$ which finishes the proof of
Proposition 7.13.

\head {\bf References}
\endhead

\enddocument